\documentclass[preprint,3p,authoryear,times]{elsarticle}

\usepackage[fleqn]{amsmath}
\usepackage[]{amsfonts}
\usepackage[]{amsthm}
\usepackage[]{amssymb}
\usepackage{empheq}
\usepackage{lscape}
\usepackage[ruled,vlined]{algorithm2e}
\usepackage{booktabs}
\usepackage{makecell}
\usepackage{tabu}
\usepackage{caption}
\usepackage{graphicx}
\usepackage[export]{adjustbox}

\usepackage[mathlines]{lineno}
\usepackage{tikz}
\usetikzlibrary{positioning, arrows.meta}
\usepackage{epstopdf}
\AppendGraphicsExtensions{.tif}

\usepackage[]{algorithm2e}
\usepackage[english]{babel}
\usepackage[T1]{fontenc}

\usepackage{hyperref} 
\usepackage{graphicx}
\usepackage{caption,subcaption}
\usepackage{tikz}
\usepackage{pstricks}
\usepackage{xstring}
\usepackage{multirow}
\usepackage{hhline}

\definecolor{darkblue}{rgb}{0,0,0.5}

\hypersetup{colorlinks=true,linkcolor=black,citecolor=darkblue,urlcolor=darkblue} 

\newcommand{\newref}[2]{\hyperref[#2]{#1~\ref*{#2}}} 

\usepackage{longtable}
\usepackage{pdflscape}
\usepackage{color,soul}
\usepackage[toc,page]{appendix}

\newcommand{\RR}{{\rm I\kern-2pt R}}
\newcommand{\Def}{\stackrel{def}{=}}

\def\Bc{{\cal B}}
\def\Cc{{\cal C}}

\def\Tc{{\cal T}}

\begin{document}

\begin{frontmatter}

\title{Dynamic System Optimum: A Projection-based Framework for Macroscopic Traffic Models}



\author[1]{Mostafa Ameli\corref{cor}} 
\cortext[cor]{Corresponding author.}
\ead{mostafa.ameli@univ-eiffel.fr}
\author[2]{Sérgio F. A. Batista}
\author[1]{Jean-Patrick Lebacque}
\author[3]{M\'onica Men\'endez}

\address[1]{COSYS, GRETTIA, University Gustave Eiffel, Paris, Île-de-France, France}

\address[2]{Department of Civil and Environmental Engineering, Imperial College London, UK}

\address[3]{Division of Engineering, New York University Abu Dhabi, Saadiyat Marina District PO Box 129188 - Abu Dhabi, United Arab Emirates}

\begin{sloppypar} 
\begin{abstract}
This paper proposes the theoretical grounds for emulating the Dynamic System Optimum with desired arrival times on regional networks, using aggregated traffic dynamics based on the Macroscopic Fundamental Diagram. We used a projection gradient-based solution method that avoids the need to compute approximations of marginal travel times. We demonstrate the application of our framework on 8-region networks and show that our approach yields improved solutions compared to classical approximation methods in the literature, such as the Method of Successive Averages and the gap-based approach.
\end{abstract}


\end{sloppypar} 
\end{frontmatter}



\begin{sloppypar} 
\section*{Highlights}
\begin{itemize}
\item We formulate the Dynamic System Optimum for the aggregated traffic models based on the Macroscopic Fundamental Diagram.
\item Our approach does not require the computation of path marginals to solve the network equilibrium. 
\item We adopt a projection-based solution algorithm to compute the global Dynamic System Optimum.
\item We showcase the application of our method on an 8-region network.
\item We show that our method provides a better solution than benchmark models, such as the MSA and gap-based approaches.
\end{itemize}

\section{Introduction}\label{Sect1}
In urban transport systems, managing congestion and minimising emissions requires an effective modelling of the dynamics of traffic and their variations over time. Dynamic Traffic Assignment (DTA) is a widely used tool in transportation planning and management that simulates how travellers make route and departure time decisions under changing network conditions. DTA models capture both individual travel behaviour and the resulting interactions within the network, reflecting demand fluctuations and the physical propagation of congestion. The theoretical foundations of traffic assignment are rooted in Wardrop’s principles \citep{Wardrop-1952}. The \textit{First Principle}, also known as \textit{User Equilibrium (UE)}, states that no traveller can unilaterally reduce travel time by changing routes. The \textit{Second Principle}, or \textit{System Optimum (SO)}, reflects a condition in which routes are allocated to minimise total travel time among all travellers in the network. In the dynamic context, analogous concepts are the Dynamic User Equilibrium (DUE) and the Dynamic System Optimum (DSO). A network reaches DUE when no traveller can reduce travel time by unilaterally switching routes at a given departure time. DSO occurs when any reallocation of route choices cannot further reduce the total travel time between all vehicles. The emulation of these equilibria in large-scale networks presents significant challenges. Solving DUE or DSO requires addressing a constrained convex optimisation problem \citep{Sheffi-1985}, where DSO further involves the computation of \textit{marginal travel times} — a computationally intensive and conceptually challenging task \citep[see e.g.][]{Ameli-2020,ameli2024collective}. In addition to the theoretical challenges, there is the added complexity of realistically replicating traffic dynamics. Accurate DTA requires time-dependent travel times that capture shockwaves, spillbacks, and dynamically changing demand, particularly in congested urban settings. These requirements pose a scalability problem for fine-grained link-based simulation \citep[see e.g.][]{Ameli-2020}.

To address the computational and scalability limitations of fine-grained simulation-based DTA models, aggregated traffic models based on the Macroscopic Fundamental Diagram (MFD) have gained prominence \citep{Daganzo-2007, Geroliminis-2008, Vickrey-2020}. These models provide an aggregated representation of traffic dynamics on a neighbourhood scale, making them suitable for large-scale applications \citep[e.g.][]{Mariotte-2020, Batista-2022}. The MFD-based framework has been successfully applied in a variety of domains, including perimeter control \citep{Haddad-2017a, Sirmatel-2019, He-2019, Haddad-2020, Sirmatel-2021, Li-2021, SIRMATEL2023104338, KUMARAGE2023104184}, pricing strategies \citep{Yang-2019, Zheng-2020, Loder-2022}, multimodal network modelling \citep{Loder-2019, Paipuri-2020a, Paipuri-2020c, BALZER2023104061}, networkwide emissions control \citep{Ingole-2020b, Batista-2022}, and ride sourcing and resilience analysis \citep{Beojone-2021, LU2024110095}.

The city or a road network is defined by a graph where the nodes represent intersections and the links of that graph represent the road segments connecting intersections. To represent the network at the macroscopic level, we define an aggregated graph where the nodes represent the regions and the links represent the allowed travel directions between adjacent regions. A fundamental component of MFD-based models is the partition of the urban network into regions where traffic is assumed to travel at the same spatial mean speed \citep{Saeedmanesh-2017, Lopez-2017, Batista-2022b}. Within each region, the traffic dynamics are governed by an intrinsic MFD that describes the relationship between the accumulation of vehicles in the region and the production of travel. This partition allows for the definition of a graph, that is, a regional graph, representing this regional network, where the nodes represent the regions and the links connecting these nodes depend on the underlying topology of the road infrastructure. Traffic is represented as exchange flows between adjacent regions. Following \cite{Batista-2021a}, we define a trip at the microscopic level as the ordered sequence of travelled links from the origin to the destination (aggregated) nodes, whereas a path in the regional graph is represented as the ordered sequence of travelled regions from the Origin (O) to the Destination (D) regions. It is important to note that many individual trips can be mapped onto a single regional path \citep{Batista-2021a}. This plays a crucial role in characterising the travel times of the paths \citep{Batista-2019b} and consequently in the calculation of DUE and DSO in regional networks. Following \cite{Batista-2019a}, paths are characterised by explicit distributions of travel distances in each region they traverse, while a trip consists of a sequential set of links with a fixed physical length. This plays an important role in the calculation of the travel times of the roads in regional networks \citep{Yildirimoglu-2014,Batista-2019b}, and therefore in the development of traffic management strategies \citep{Yildirimoglu-2014,Ramezani-2015,Yildirimoglu-2018,Batista-2019b}. 

Several computational frameworks have been proposed to compute DUE on regional networks under MFD dynamics. Early works such as \citet{Yildirimoglu-2014} and \citet{Batista-2019b} rely on the Method of Successive Averages (MSA) or Monte Carlo simulation, incorporating distance distributions and regional speed evolution. Later studies expanded this to consider bounded rationality \citep{Batista-2020b} and accumulation-based constraints \citep{HUANG20201}. From a more analytical perspective, \citet{Leclercq-2013} and \citet{Laval-2017} derived closed-form equilibrium solutions for simplified two-region and a highway network. \citet{Yildirimoglu-2015} discussed a perimeter control framework with route guidance that targets the DSO conditions of the regional network. The authors determined the network equilibrium using the MSA and considered the calculation of the marginal travel times. From a different perspective, \cite{AGHAMOHAMMADI2020101} introduced a continuum-space DSO framework, aligned with the MFD traffic flow theory. \citet{Aghamohammadi-2020} provides a comprehensive review of the literature on these continuum models.

From a different perspective, several studies have addressed the dynamic choice modelling of departure time within the MFD framework. \citet{ARNOTT2018150} proposed a two-loop method to solve the equilibrium in the bathtub model, considering entry rates and commuter populations. \citet{LAMOTTE2018794} showed that a first-in first-out pattern emerges under inelastic demand and heterogeneous travel distances, while \citet{LAMOTTE202114} found that such equilibria may be unstable with scheduling preferences. \citet{Amirgholy-2017} integrated Vickrey’s model with MFD dynamics to optimise the design of the transit and minimise travel, operational and environmental costs. Using mean-field games, \citet{ameli2022departure} solved the departure-time user equilibrium in accumulation-based MFDs. \citet{YILDIRIMOGLU2021103391} proposed staggering work schedules through multi-objective optimisation to alleviate congestion and delay. \citet{ZHONG2021103190} derived departure-time DUE in isotropic trip-based MFDs, showing that inflow constraints can prevent hypercongestion. \citet{Zhong9102445} studied DSO with joint modelling of departure time and route choices across multiregional MFD networks, which avoids strong assumptions such as linearisation of the MFD.

The development of simulation-based dynamic traffic assignment models for regional networks is important not only for the design of traffic management strategies, but also for the setting of environmental policies, to name a few examples, as evidenced in the review work done by \citep{JOHARI2021103334}. Despite the extensive literature on dynamic system optimum in both microscopic and regional traffic models, no existing work, to our knowledge, has provided a theoretical foundation for emulating simulation-based DSO on regional MFD networks without relying on approximate marginal travel times. For this, we resort to a projected gradient solution algorithm to solve the DSO, which avoids the calculation of explicit marginal times. The framework proposed in this paper can handle time-varying demand profiles, desired arrival times, and region-to-region path flows within a dynamic simulation. We showcase the application of this framework on a 4- and 8-region network and analyse the DSO conditions against benchmark approximations in the literature. 

This paper is organised as follows. Section \ref{Sect2} introduces the regional MFD dynamics with the desired arrival times and the mathematical notation. Section \ref{Sect3} formulates the mathematical formulation of the Dynamic System Optimum on regional networks considering the MFD dynamics and the desired arrival times. Section \ref{Sect4} describes the solution algorithm followed to solve for the DSO conditions in regional networks. Section \ref{Sect5} introduces the projection method. Section \ref{Sect6} presents the application of this theoretical framework on 4-region and 8-region networks and analyses the results against the benchmark approximations used in the literature. Section \ref{Sect7} summarises the main conclusions of this paper and discusses potential future research lines. 

\section{Regional MFD-based dynamics with desired arrival times}\label{Sect2}

\subsection{Regional networks and demand}\label{Sect2.1}
The graph $\mathcal{G}=(R,U)$ defines the regional network where the nodes $r \in R$ represent the regions and $u \in U$ is the set of edges indicating the travel directions. Thus, an edge $u=(i,j)\in U$ corresponds to two regions $(i)$ and $(j)$ that are neighbours such that travellers can travel from $(i)$ to $(j)$. In this situation, we denote $j \in \Gamma^+(i)$ and $i \in \Gamma^-(j)$). The cardinality, that is, the number of elements in each set of $R$ and $U$ is $\bar{R}$ and $\bar{U}$, respectively.

A path in the regional network is defined as an ordered sequence of travelled regions between the origin (O) and the destination (D) \citep{Batista-2021a}. Any region of the regional network can be an Origin or Destination for drivers, but also an intermediate region when drivers travel through it from one preceding region to the next region. We define $\mathcal{O}$, $\mathcal{D}$, and $\mathcal{J}$ as the sets of Origin, Destination, and Intermediate regions in the regional network, respectively. We distinguish origin and destination regions from intermediate regions as they require different management of the inflow and outflow of vehicles. Drivers also travel through intermediate regions when travelling from a preceding region to the next region to be travelled in the sequence of the path. We define $\Gamma^+(i)$ and $\Gamma^-(i)$ as the sets of successor and predecessor regions that travelled before region $i$ and that would be travelled after region $i$, respectively.

We define $\Delta^{OD,t_a}(t)$ as the total disaggregated demand at time $t \in [0,T]$ travelling between regions O and D, and with a desired arrival time $t_a \in \mathcal{T}_a$. The term $\mathcal{T}_a$ represents a finite set of desired arrival times, and $T$ denotes the total simulation period. The total demand $\Delta^{OD,t_a}$ is:

\begin{equation}\label{Eq_TotalDemand1}
\Delta^{OD,t_a} = \int_{\mathcal{T}_d} \Delta^{OD,t_a}(t) \,\mathrm{dt}, \quad \forall O,D \wedge \forall t_a \in\Tc_a
\end{equation}

\noindent where $\Delta^{OD,t_a}(.) \ge 0$; and $\mathcal{T}_d$ represents the set of departure times between $[0,T]$, with $T$ the total simulation period.

\subsection{Description of the inflows and outflows with MFD dynamics}\label{Sect2.2}
This section describes the calculation of inflows and outflows in regions using the MFD dynamics. Note that the inflows in a given region $i$ are the outflows of the predecessor regions $j \in \Gamma^-(i)$, and vice versa. Therefore, our implementation of the MFD dynamic model consists of three loops. The first loop calculates demand and supplies. The second loop calculates the inflows (and thus outflows) in the regions. The third loop determines the evolution of traffic dynamics in the regions. These three loops constitute the state equations of the system. 

\subsubsection{Supplies and demands}\label{Sect2.2.1}
We distinguish between the origin, destination, and intermediate regions. We assume the existence of an MFD for each region $i \in \mathcal{J} \cup \mathcal{O}$. If $i \in \mathcal{J}$, it is an effective MFD calculated based on traffic data. For $O \in \mathcal{O}$, the MFD represents the function of $O$ as a buffer between the demand $\Delta^{OD,t_a}$ and the supply of the regions $i \in \Gamma^+(O)$. 

\begin{figure}[h!]
    \centering
    \subfloat[\centering]{{\includegraphics[width=0.35\textwidth]{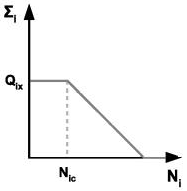}}}%
    \qquad
    \subfloat[\centering]{{\includegraphics[width=0.35\textwidth]{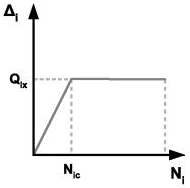}}}
    \caption{(a) Demand as a function of $N_i$. (b) Supply as a function of $N_i$.}
    \label{Fig:demandandsupply}
\end{figure}

For a generic region $i \in \mathcal{J}$, the MFD results in the supply $\Sigma_i(.)$ and demand $\Delta_i(.)$ functions. \hyperref[Fig:demandandsupply]{Figure~\ref*{Fig:demandandsupply}} shows the demand and supply as a function of the total number of travellers $N_i$ (accumulation).

There are two state variables. Let $N_i^{D,t_a}(t)$ be the number of travellers in the region $i$ at time $t$ with destination $D$ and desired arrival time $t_a$. Then, the total number of travellers $N_i(t)$ present in the region $i$ at the time instant $t$ is given by the following:

\begin{equation}\label{Eq_Ni}
N_i(t) = \sum_{\substack{D \in \mathcal{D} \\ t_a \in \mathcal{T}_a}} N_i^{D,t_a}(t), \forall t
\end{equation}

The total supply $\sigma_i(t)$ and demand $\delta_i(t)$ for the region $i$ are determined as follows:

\begin{equation}\label{Eq_Demand1}
\sigma_i(t) = \Sigma_i(N_i(t)), \forall i \in \mathcal{J}
\end{equation}

\begin{equation}\label{Eq_Supply1}
\delta_i(t) = \Delta_i(N_i(t)), \forall i \in \mathcal{J}
\end{equation}

\begin{figure}[h!]
    \centering
    \includegraphics[width=0.5\textwidth]{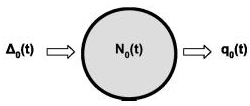}
    \caption{Dynamics of a buffer zone.}
    \label{Fig:demandandsupply_1}
\end{figure}

\begin{figure}[h!]\centering
    \includegraphics[width=0.35\textwidth]{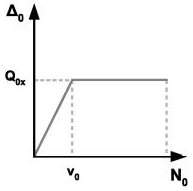}
    \caption{Demand function $\Delta_O (N_O)$ of the buffer.}
    \label{Fig:demandandsupply_2}
\end{figure}

Recall that the origin regions $O \in \mathcal{O}$, which are buffers, represent a special case. Their functioning can be summarised as follows. We also define two state variables. Let $N^{OD,t_a}(t)$ be the number of travellers in the buffer zone $O \in \mathcal{O}$ at time instant $t$, with the destination region $D$ and the desired arrival time $t_a$. \hyperref[Fig:demandandsupply_1]{Figure~\ref*{Fig:demandandsupply_1}} depicts a schematic representation of a buffer origin region, with input demand $\Delta^{OD,t_a}(t)$ and output $q^O(t)$. \hyperref[Fig:demandandsupply_2]{Figure~\ref*{Fig:demandandsupply_2}} shows the demand profile $\Delta_0$.

\noindent The total number of travellers $N^O(t)$, at time $t$, should satisfy the following condition:

\begin{equation}
\dot{N}^O(t) = \begin{cases} \Delta^O(t) - q^O(t), & \text{if } N^O(t) > 0 \\
\text{max}\left(0,\Delta^O (t) - q^O(t)\right), & \text{if } N^O(t) = 0
\end{cases}
\end{equation}

\noindent where $q^O(t)$ equals the maximum flow of buffer $O$ (if $N^O(t) \not= 0$) or $\Delta^O (t)$ (if $N^O(t) = 0$). Alternatively, the outflow $q^O(t)$ can be calculated as a function of the downstream supply of $O$ and the demand $\delta^O (t)$ of the buffer. This approach enables us to take into account multiple destinations. So, in order to evaluate $\delta_O(t)$ we propose the following demand function $\Delta_O(N^O)$ for the buffer:

\begin{equation}
\Delta_O(N^O(t)) = \begin{cases}
Q_{Ox}, & \text{if } N^O(t) \ge \nu_O \\
Q_{Ox} \cdot \frac{N^O(t)}{\nu_O}, & \text{if } N^O(t) < \nu_O \\
\end{cases}
\end{equation}

\noindent where $Q_{Ox}$ should be chosen large. For example, $Q_{Ox}$ should be the sum of the maximum inflows of $i \in \Gamma^+(i)$. Then, we define the demand of the buffer zone, $\delta_O(t)$ as follows:

\begin{equation}\label{Eq_SupplyO}
\delta_O(t) = \Delta_O(N^O(t)), \forall O \in \mathcal{O} \wedge \forall t  
\end{equation}

Note that this alternative option entails a systematic error in a non-empty buffer with a non-zero travel time. Typically, the error would be of size $\nu_O$, and the error in travel time would be approximately $\nu_O/Q_{Ox}$ per traveller (for an upper bound).

Finally, for the special case of destination regions $D \in \mathcal{D}$, we require as data the exit supply $\sigma_D$.

\subsubsection{Disaggregation of supplies and demands: the calculation of path flows}\label{Sect2.2.2}
We start by calculating the path flows for a generic region $i \in \mathcal{J}$. This means calculating the inflows from a generic predecessor region $j$ to region $i$ for $i \in \mathcal{J}$ and $j \in \Gamma^-(i)$. In the first step, we calculate the partial demands $\delta_{i\ell }, \forall \ell  \in \Gamma^+(i) \wedge i \in \mathcal{J}$. The number of drivers with final destination $D \in \mathcal{D}$ and arrival time $t_a \in \mathcal{T}_a$, who wish to travel from a generic region $i$ to a next adjacent region $\ell $, is given by $\gamma_{i\ell }^{D,t_a}(t) \cdot N_i^{D,t_a}(t)$ at time $t$. The term $\gamma_{i\ell }^{D,t_a}(t)$ represents the fraction of travellers in the region $i$ who travel to $\ell $ as their next region at time $t$, given their destination region $D$ and the desired arrival time $t_a$. Therefore, the proportion of travellers in the region $i$ who want to go to $\ell  \in \Gamma^+(i)$ at time $t$ is $\sum_{\substack{D \in \mathcal{D}, t_a \in \mathcal{T}_a}} \gamma_{i\ell }^{D,t_a}(t)\cdot N_i^{D,t_a}(t)/N_i(t)$. Based on this, we define the demand in the region $i$, $\delta_{i\ell }(t)$, as follows: 

\begin{equation}\label{Eq_Supplyij}
\delta_{i\ell }(t) = \delta_i(t) \cdot \frac{\sum_{\substack{D \in \mathcal{D}, t_a \in \mathcal{T}_a}} \gamma_{i\ell }^{D,t_a}(t)\cdot N_i^{D,t_a}(t)}{N_i(t)}, \forall i \in \mathcal{J} \cup \mathcal{O} \wedge \ell \in \Gamma^+(i) 
\end{equation}

Given this, we use the STRADA model \citep{buisson1995macroscopic} to calculate partial supplies and flows. This model assumes that the supply for traffic entering the region $j$ and coming from the predecessor region $i$, $\sigma_{ij}$, is a fraction of the total supply of the region $j$:

\begin{equation}\label{Eq_DemandSTRADA}
\sigma_{ij} = \beta_{ij} \cdot \sigma_j, \forall i \in \Gamma^-(j) \wedge j \in \mathcal{J} \cup \mathcal{D}  
\end{equation}

\noindent with constant split coefficients $\beta_{ij}$.

This phenomenological model for complex intersections does not satisfy the invariance principle \citep{lebacque2005first}. The inflow of the region $j$ coming from the predecessor region $i$, $q_{ij}(t)$, is given by the minimum value between the supply $\sigma_{ij}(t)$ and the demand $\delta_{ij}(t)$:

\begin{equation}\label{Eq_Inflows1}
q_{ij}(t) = \mathrm{min} \left(\sigma_{ij}(t),\delta_{ij}(t)\right), \quad \forall i \in \Gamma^-(j) \wedge j \in \mathcal{J} \cup \mathcal{D}
\end{equation}

\subsubsection{Optimization model for flows}\label{Sect2.2.3}
An alternative model for region inflows which satisfies the invariance principle is the flow optimisation model \citep{lebacque2005first}. This model does not propose the calculation of partial supplies, but proposes a direct calculation of flows $q_{ji}, \forall j \in \Gamma^-(i) \wedge \forall i \in \mathcal{J} \cup \mathcal{D}$.

Calculating inflows into the region $i$ is achieved by maximising the objective function $Z$:

\begin{equation}\label{Eq_Optim1_KKT}
\textrm{max } Z = \sum\limits_{\ell \in \Gamma^-(i)} \Phi_{\ell i}(q_{\ell i})    
\end{equation}

\noindent subject to:

\begin{equation}\label{Eq_Optim2_KKT}
0 \le q_{\ell i} \le \delta_{\ell i}(t), \forall \ell  \in \Gamma_i^-  
\end{equation}

\begin{equation}\label{Eq_Optim3_KKT}
\sum\limits_{\ell  \in \Gamma_i^-} q_{\ell i} \le \sigma_i(t)
\end{equation}

As usual for $\Phi_{\ell i}(q_{\ell i})$, we choose a concave quadratic increasing function:

\begin{equation}\label{Eq_Optim4_KKT}
\Phi_{\ell i}(q_{\ell i}) = q_{\ell i} - \frac{q_{\ell i}^2}{q_{\ell i,x}}, \forall \ell  \in \Gamma^-(i) 
\end{equation}

\noindent where $q_{\ell i,x}$ is the maximum flow from region $\ell $ to $i$.

The idea of this model can be traced back to \cite{holden1995mathematical}; in the case of a simple merge, the model can be shown to be equivalent to the merge model in \cite{daganzo1995cell}.

To solve the optimisation problem defined in Equations \hyperref[Eq_Optim1_KKT]{\ref*{Eq_Optim1_KKT}} to \hyperref[Eq_Optim4_KKT]{\ref*{Eq_Optim4_KKT}}, we write the optimality conditions of Karush-Kuhn-Tucker (KKT). Then the gradient KKT condition yields the following:

\begin{equation}\label{Eq_KKT1}
1 - \frac{q_{\ell i}}{q_{\ell i,x}} = \zeta_i + \varphi_{i\ell } - \psi_{i\ell }
\end{equation}

\noindent $\zeta_i$ is the KKT coefficient of the demand constraint $\sigma_i$ (i.e., Equation \hyperref[Eq_Optim3_KKT]{\ref*{Eq_Optim3_KKT}}); $\varphi_{i\ell }$ is the KKT coefficient of the $q_{i\ell } \le \delta_{i\ell }$ constraint (i.e. Equation \hyperref[Eq_Optim2_KKT]{\ref*{Eq_Optim2_KKT}}); and $\psi_{i\ell }$ is the KKT coefficient of the $0 \le q_{i\ell}$ constraint (i.e. Equation \hyperref[Eq_Optim2_KKT]{\ref*{Eq_Optim2_KKT}}).

The KKT conditions include Equation \hyperref[Eq_KKT1]{\ref*{Eq_KKT1}}, the constraints Equations \hyperref[Eq_Optim2_KKT]{\ref*{Eq_Optim2_KKT}} to \hyperref[Eq_Optim4_KKT]{\ref*{Eq_Optim4_KKT}}, positivity constraints and complementarity constraints for the KKT coefficients $\phi_{i\ell}$ and $\psi_{i\ell }$. Given three real numbers $a,b,x$ with $a\leq b$, we define $P_{[a,b]}(x)$ as $a$ if $x\leq a$, $x$ if $a\leq x \leq b$ and $b$ if $x\geq b$. With this definition, we solve the KKT optimality conditions for $\phi_{i\ell}$ and $\psi_{i\ell}$, resulting in the following:

\begin{equation}\label{Eq_KKT2}
q_{\ell i}(t) = P_{[0,\delta_{\ell i}(t)]} \cdot \left(q_{\ell i,x}\cdot (1 - \zeta_i(t))\right)    
\end{equation}

\noindent and we solve $\zeta_i(t)$ by:

\begin{equation}\label{Eq_KKT3}
\sum_{\ell  \in \Gamma^-(i)} q_{\ell i}(t) = \textrm{min}\Bigg(\sigma_i(t),\sum_{\ell  \in \Gamma^-(i)} \delta_{\ell i}(t)\Bigg)    
\end{equation}

In equation \hyperref[Eq_KKT3]{\ref*{Eq_KKT3}} we choose the term $\textrm{min}\Big(\sigma_i(t),\sum_{k \in \Gamma^-(i)} \delta_{\ell i}(t)\Big)$, such that the solution $\zeta_i(t)$ is unique and continuously depends as a Lipschitz piecewise affine function on $\sigma_i(t)$ (and on $\delta_{\ell i}(t)$). To guarantee unicity, we complete the system of equations \hyperref[Eq_KKT2]{\ref*{Eq_KKT2}} to \hyperref[Eq_KKT3]{\ref*{Eq_KKT3}} with:

\begin{equation}\label{Eq_KKT4}
\zeta_i(t) \in \Bigg[1-\max\limits_{k \in \Gamma^-(i)}\Bigg(\frac{\delta_{\ell i}(t)}{q_{\ell i,x}} \Bigg),1 \Bigg]
\end{equation}

\begin{figure}[h!]\centering
    \includegraphics[width=0.45\textwidth]{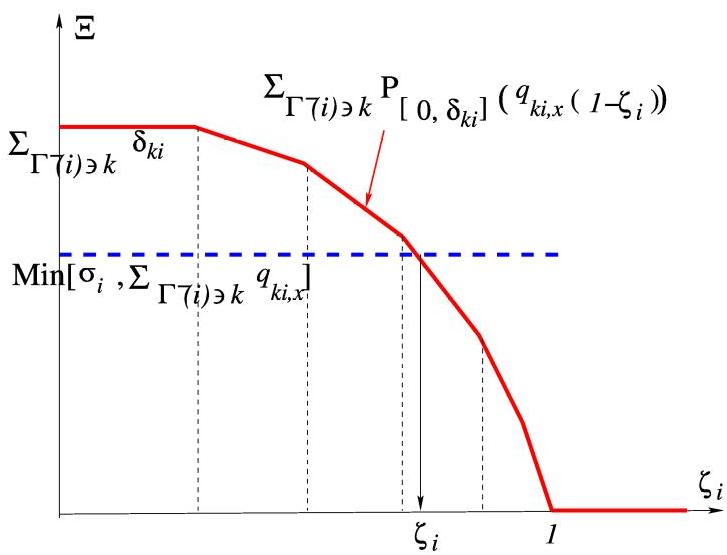}
    \caption{Intersection of $\Xi(\zeta_i)$ with the line $\textrm{min}\Bigg(\sigma_i(t),\sum_{\ell  \in \Gamma^-(i)} \delta_{\ell i}(t)\Bigg)$.}
    \label{Fig:KKT_2}
\end{figure}

\noindent The system of equations \hyperref[Eq_KKT2]{\ref*{Eq_KKT2}} to \hyperref[Eq_KKT4]{\ref*{Eq_KKT4}} must be solved with respect to $\zeta_i(t)$. The following function, \hyperref[Eq_Xi]{\ref*{Eq_Xi}} is decreasing and linear piecewise. It is illustrated in Figure \hyperref[Fig:KKT_2]{\ref*{Fig:KKT_2}} and is strictly decreasing in the interval $\Bigg[1-\max\limits_{\ell  \in \Gamma^-(i)}\Bigg(\frac{\delta_{\ell i}(t)}{q_{\ell i,x}} \Bigg),1 \Bigg]$. Thus, its intersection with the line $\textrm{min}\Bigg(\sigma_i(t),\sum_{\ell  \in \Gamma^-(i)} \delta_{\ell i}(t)\Bigg)$ is unique.

\begin{equation}\label{Eq_Xi}
    \Xi(\zeta_i)\Def  \zeta_i \rightarrow  \sum_{k \in \Gamma^-(i)} P_{[0,\delta_{\ell i}(t)]} \cdot \left(q_{\ell i,x}\cdot (1 - \zeta_i  )\right)  
\end{equation}

We solve this system of equations by dichotomy since the interval in which the solution lies is known, or by using the piecewise linearity of the function $\Xi(\zeta_i): \zeta_i \rightarrow  \sum_{\ell  \in \Gamma^-(i)} P_{[0,\delta_{\ell i}(t)]} \cdot \left(q_{\ell i,x}\cdot (1 - \zeta_i  ) \right)$. The latter approach involves making the list of cut values, that is, the values for which the slope of the function $\zeta_i \rightarrow  \sum_{\ell  \in \Gamma^-(i)} P_{[0,\delta_{\ell i}(t)]} \cdot (q_{\ell i,x}\cdot (1 - \zeta_i  )$ changes. Then we sort the cut values in increasing order and calculate the values of $\sum_{\ell  \in \Gamma^-(i)} P_{[0,\delta_{\ell i}(t)]} \cdot (q_{\ell i,x}\cdot (1 - \zeta_i(t)))$ for the cut values. These values are decreasing and can be compared to $\textrm{min}\Bigg(\sigma_i(t),\sum_{\ell \in \Gamma^-(i)} \delta_{\ell i}(t)\Bigg)$. We deduce in which interval of the cut values lies the solution. In this interval, the system of equations \hyperref[Eq_KKT2]{\ref*{Eq_KKT2}} to \hyperref[Eq_KKT4]{\ref*{Eq_KKT4}} is linear and is solved trivially. This procedure produces an exact solution in a finite number of steps.


\subsubsection{Calculation of the traffic states}\label{Sect2.2.4}
This section describes the calculation of the dynamics of the number of vehicles for all regions $i \in \mathcal{O} \cup \mathcal{J}$, and the special cases of the Origin $O \in \mathcal{O}$ and Destination $D \in \mathcal{D}$ regions.

We start with the regions $i \in \mathcal{J}$. The number of vehicles travelling from the region $i$ to the next adjacent region $j\in \Gamma^+(i)$ at time $t$ and that have a destination $D \in \mathcal{D}$ with a desired arrival time $t_a \in \mathcal{T}_a$, $N_{ij}^{D,t_a}(t)$, is given by 
$$ \gamma_{ij}^{D,t_a}(t)\cdot N_{i}^{D,t_a}(t) $$

Thus the fraction of travelers in region $i$ bound for adjacent region $j\in \Gamma^+(i)$ at time $t$, which have destination $D \in \mathcal{D}$ and desired arrival time $t_a \in \mathcal{T}_a$ is calculated as follows: 
\begin{equation}\label{Eq_qijdta}
 \frac{\gamma_{ij}^{D,t_a}(t)\cdot N_{i}^{D,t_a}(t)}{\sum\limits_{\delta\in \mathcal{D}, \theta\in \mathcal{T}_a} \gamma_{ij}^{\delta,\theta}(t)\cdot N_{i}^{\delta,\theta}(t)}, \quad \forall   j\in \Gamma^+(i)\wedge D \in \mathcal{D}\wedge t_a \in \mathcal{T}_a
\end{equation}

We apply this fraction to the flow $q_{ij}(t)$ of vehicles travelling from regions $i$ to $j$ to disaggregate $q_{ij}(t)$ regarding the destination $D \in \mathcal{D}$ and the desired arrival time $t_a \in \mathcal{T}_a$ is given by:

\begin{equation}\label{Eq_qijDta}
q_{ij}^{D,t_a}(t) = q_{ij}(t) \cdot \frac{\gamma_{ij}^{D,t_a}(t)\cdot N_{i}^{D,t_a}(t)}{\sum\limits_{\delta\in \mathcal{D}, \theta\in \mathcal{T}_a} \gamma_{ij}^{\delta,\theta}(t)\cdot N_{i}^{\delta,\theta}(t)} , \quad \forall D \in \mathcal{D} \wedge \forall t_a \in \mathcal{T}_a   
\end{equation}

The evolution of the number of vehicles $N_{i}^{D,t_a}(t)$ and $N_{i}(t)$ as a function of time dictates the traffic dynamics in the region $i$, which is given by the following system of first-order differential equations:

\begin{equation}\label{Eq_DNiDta}
\dot{N}_{i}^{D,t_a}(t) = \sum_{\ell  \in \Gamma^-(i)} q_{\ell i}^{D,t_a}(t) - \sum_{j \in \Gamma^+(i)} q_{ij}^{D,t_a}    
\end{equation}

\begin{equation}\label{Eq_DNi}
\dot{N}_{i}(t) = \sum_{\ell  \in \Gamma^-(i) } q_{\ell i}(t) - \sum_{j \in \Gamma^+(i)} q_{ij}
\end{equation}

Note that Equation \hyperref[Eq_DNi]{\ref*{Eq_DNi}} is obtained by summing Equation \hyperref[Eq_DNiDta]{\ref*{Eq_DNiDta}} over all $D \in \mathcal{D}$ and $t_a \in \mathcal{T}_a$. This is also consistent with Equation \hyperref[Eq_Ni]{\ref*{Eq_Ni}}.

We focus on the special case of the Origin regions $O \in \mathcal{O}$. Similarly to the previous case, from Equation \hyperref[Eq_qijDta]{\ref*{Eq_qijDta}}, the flow of vehicles travelling from the origin $O \in \mathcal{O}$ region to the next adjacent region $i$ is the following:

\begin{equation}
q_{Oi}^{D,t_a}(t) = q_{Oi}(t)\cdot  \frac{\gamma_{Oi}^{D,t_a}(t)\cdot N_{O}^{D,t_a}(t)}{\sum\limits_{\delta\in \mathcal{D}, \theta\in \mathcal{T}_a } \gamma_{Oi}^{\delta,\theta}(t)\cdot N_{O}^{\delta,\theta}(t)}, \quad \forall  i \in\Gamma^+(O) \wedge \forall  D \in \mathcal{D}  \wedge \forall t_a \in \mathcal{T}_a    
\end{equation}

In addition, Equation \hyperref[Eq_DNi]{\ref*{Eq_DNi}} must be satisfied with the buffer input being the upstream demand. That is,

\begin{equation}\label{Eq_DNODta}
\dot{N}_{O}^{D,t_a}(t) = \Delta_O^{D,t_a}(t) - \sum_{i \in \Gamma^+(O)} q_{Oi}^{D,t_a}, \forall O \in \mathcal{O} \wedge \forall D \in \mathcal{D} \wedge \forall t_a \in \mathcal{T}_a    
\end{equation}

\begin{equation}\label{Eq_DNO}
\dot{N}_{O}(t) = \Delta_O(t) - \sum_{i \in \Gamma^+(O)} q_{Oi} (t)
\end{equation}

\noindent where $\Delta_O(t) = \sum_{\substack{\delta \in {\mathcal{D} } \\ \theta \in \mathcal{T}_a}} \Delta_O^{\delta,\theta}(t)$. Similarly to the previous case, Equation \hyperref[Eq_DNO]{\ref*{Eq_DNO}} is obtained by summing Equation \hyperref[Eq_DNODta]{\ref*{Eq_DNODta}} in all $D \in \mathcal{D}$ and $t_a \in \mathcal{T}_a$, which is consistent with Equation \hyperref[Eq_Ni]{\ref*{Eq_Ni}}.

\section{System Optimum on regional networks with desired arrival time}\label{Sect3}
This section formulates the optimal conditions of the dynamic system in regional networks, where travellers have a desired arrival time.

\subsection{Calculation of the System Optimum}\label{Sect3.1}
The dynamical system includes the state variables and subsidiaries, the decision variables, the state equations, and the constraints. The state variables are $N(t): N_i^{D,t_a}(t), \forall i \in \mathcal{J} \cup \mathcal{O} \wedge D \in \mathcal{D} \wedge t_a \in \mathcal{T}_a$. The subsidiary variables are $N_i(t), \forall i \in \mathcal{J} \cup \mathcal{O}$. The decision variables are $\Delta(t): \Delta_O^{D,t_a}, \forall O \in \mathcal{O} \wedge \forall D \in \mathcal{D} \wedge \forall t_a \in \mathcal{T}_a$ and $\gamma(t): \gamma_{ij}^{D,t_a}, \forall i \in \mathcal{J} \cup \mathcal{O} \wedge \forall j \in \Gamma^+(i) \wedge \forall D \in \mathcal{D} \wedge \forall t_a \in \mathcal{T}_a$.

The state equations of the system are described by the following: Equations \hyperref[Eq_Demand1]{\ref*{Eq_Demand1}}, \hyperref[Eq_Supply1]{\ref*{Eq_Supply1}}, \hyperref[Eq_SupplyO]{\ref*{Eq_SupplyO}} and \hyperref[Eq_Supplyij]{\ref*{Eq_Supplyij}} represent the supplies and demands of the system; Equations \hyperref[Eq_DemandSTRADA]{\ref*{Eq_DemandSTRADA}} and \hyperref[Eq_Inflows1]{\ref*{Eq_Inflows1}}  or the system of Equations \hyperref[Eq_KKT2]{\ref*{Eq_KKT2}}, \hyperref[Eq_KKT3]{\ref*{Eq_KKT3}} and \hyperref[Eq_KKT4]{\ref*{Eq_KKT4}} describe the inflow dynamics in the regions; Equations \hyperref[Eq_qijDta]{\ref*{Eq_qijDta}}, \hyperref[Eq_DNiDta]{\ref*{Eq_DNiDta}}, \hyperref[Eq_DNi]{\ref*{Eq_DNi}}, \hyperref[Eq_DNODta]{\ref*{Eq_DNODta}} and \hyperref[Eq_DNO]{\ref*{Eq_DNO}} describe the evolution of the accumulation of vehicles, i.e. the traffic dynamics, as a function of time in the network; and Equation \hyperref[Eq_Ni]{\ref*{Eq_Ni}} represents the total accumulation in the network as a function of time.

This dynamical system is subject to various constraints. The constraints applied to the state variables are built into the model by the effect of the supply and demand functions.  The constraints applied to the command variables concern both the assignment coefficients $\gamma_{ij}^{D,t_a}$ and the total demand $\Delta_O^{D,t_a}(t)$. The constraints on $\gamma_{ij}^{D,t_a}$ are:

\begin{equation}\label{Eq_Constraint1}
\sum_{j \in \Gamma^+(i) } \gamma_{ij}^{D,t_a}(t) = 1, \quad \forall i \in \mathcal{J} \cup \mathcal{O} \wedge \forall j \in \Gamma^+(i) \wedge \forall D \in \mathcal{D} \wedge \forall t_a \in \mathcal{T}_a \wedge t \in [0,T]      
\end{equation}

\begin{equation}\label{Eq_Constraint2}
\gamma_{ij}^{D,t_a}(t) \ge 0, \quad \forall i \in \mathcal{J} \cup \mathcal{O} \wedge \forall j \in \Gamma^+(i) \wedge \forall D \in \mathcal{D} \wedge \forall t_a \in \mathcal{T}_a \wedge t \in [0,T]
\end{equation}
which express that for all $t,t_a,D,i$ the coefficients $\gamma_{ij}^{D,t_a}(t)$ can be viewed as the probability of chosing $j$.

We also have the constraint on the total demand given by Equation \hyperref[Eq_TotalDemand1]{\ref*{Eq_TotalDemand1}}. These constraints bound the decision variables of the optimisation problem. We denote:

\begin{equation}\label{Eq_Constraint3}
\Delta \in \mathcal{B} \,  \mbox{ with   }  \mathcal{B} = \times_{O\in\mathcal{O},D\in\mathcal{D},t_a\in\Tc_a} \mathcal{B}_O^{D,t_a} 
\end{equation}
and for all $O,D,t_a$, 
\begin{equation*}
    \mathcal{B}_O^{D,t_a} = \left\{ \Delta_O^{D,t_a} (.) \; / \, \int_0^T dt\, \Delta  _O^{D,t_a} (t) = \overline{ \Delta  }_O^{D,t_a} \vee  \Delta_O^{D,t_a} (.) \geq 0 \right\} 
\end{equation*}

\begin{equation}\label{Eq_Constraint4}
\gamma (t) \in \Cc \quad \forall t \in [O,t_f]  
\end{equation}

Inflow is defined within the time interval $[0,T]$. Therefore, $\Delta_O^{D,t_a}(t)=0, \forall t \ge T$. However, the network is still not empty for $t\ge T$. If the network inflow is nil after $T$, the network will be empty in finite time, say after $t_f$ (final time). We might estimate $t_f$. Alternatively, we might fix it $t_f>T$ and penalise travellers who are still inside the network at $t_f$.

All in all, we can summarise the state equations of the system as follows:

\begin{equation}\label{Eq_Optim1}
\dot{N}(t) = F(N(t),\Delta(t),\gamma(t))    
\end{equation}

\noindent subject to the following set of constraints:





\begin{align}\label{Eq_Combined}
&\Delta \in \mathcal{B}   \nonumber \\
&\gamma(t) \in \mathcal{C}, &&\forall t \in [0,t_f] \\
&\Delta(t)=0, &&\forall t \ge T \nonumber \\
&N(0) = N_{\text{init}} \nonumber
\end{align}

\noindent where $N_{init}$ is the initial number of vehicles already present in the network, which may be empty. Demand $\Delta(t)=0, \forall t \ge T$ might be replaced by additional data, for example, low traffic assignment or similar out-of-peak data. 

\subsection{Criteria for System Optimum}\label{Sect3.2}
Considering the total time spent ($TTS$) in the regional network, we can write:

\begin{equation}\label{Eq_TTS1}
TTS = \sum_{i \in \mathcal{J} \cup \mathcal{O}} \int_{0}^{t_f} N_i(t)\cdot \mathrm{dt}  
\end{equation}

We can also introduce the late/early arrival cost/penalty. It is necessary to consider a destination $D \in \mathcal{D}$, with a flow $q_{iD}^{D,t_a}(t)$ of travellers arriving through the region $i \in \Gamma^-(D)$. Let $\mathcal{L}(t_a,t)$ be the cost incurred by a traveller who arrives at time $t$ where her/his desired arrival time is $t_a$. The total cost incurred for early/late arrival time is then:

\begin{equation}
TAC = \sum_{\substack{D \in {\mathcal{D}} \\ i \in \Gamma^-(D)}} q_{iD}^{D,t_a} \mathcal{L}(t_a,t) \mathrm{dt}
\end{equation}

\noindent where $TAC$ is the cost of late arrival; $\mathrm{dt}$ is the time step. So, the total early/late arrival cost is: 

\begin{equation}\label{Eq_TAC}
TAC = \sum\limits_{\substack{D \in {\mathcal{D} } \\ i \in \Gamma^-(D)}} \sum_{t_a \in \mathcal{T}_a} \int_0^{t_f} q_{iD}^{D,t_a} \mathcal{L}(t_a,t) \mathrm{dt}      
\end{equation}

\noindent where $\mathcal{L}(t_a,t)$ could be a piecewise linear function of $t-t_a$ as used for the bathtub dynamic traffic assignment problem \citep{vickrey1969congestion,fosgerau2015congestion,ameli2022departure}, or a convex piecewise differentiable function of $t-t_a$.

Note that $q(t) \stackrel{\text{def}}{=} q_{iD}^{D,t_a}, \forall D \in \mathcal{D} \wedge i \in \Gamma^-(D) \wedge t_a \in \mathcal{T}_a$ is a vector that is a function of $N(t)$ through Equations \hyperref[Eq_Demand1]{\ref*{Eq_Demand1}}, \hyperref[Eq_Supply1]{\ref*{Eq_Supply1}}, \hyperref[Eq_SupplyO]{\ref*{Eq_SupplyO}}, \hyperref[Eq_Supplyij]{\ref*{Eq_Supplyij}}, and either the STRADA Equations \hyperref[Eq_DemandSTRADA]{\ref*{Eq_DemandSTRADA}} and \hyperref[Eq_Inflows1]{\ref*{Eq_Inflows1}}  or the system of Equations \hyperref[Eq_KKT2]{\ref*{Eq_KKT2}}, \hyperref[Eq_KKT3]{\ref*{Eq_KKT3}} and \hyperref[Eq_KKT4]{\ref*{Eq_KKT4}}, that is:

\begin{equation}
q(t) = Q(N(t),\gamma(t))    
\end{equation}

Finally, we introduce a terminal criterion $TC$ to ensure that the network is not too full at the end of the period $[0,t_f]$, that is:

\begin{equation}\label{Eq_TC}
TC \stackrel{\text{def}}{=} \sum\limits_{i \in \mathcal{J} \cup \mathcal{O}} \mathcal{M}_i\left(N_i(t_f)\right)     
\end{equation}

We could choose $\mathcal{M}_i(N_i(t_f)) = \frac{1}{2} \cdot M_i \cdot N_i^2 (t_f)$. Thus, the criterion is the sum of the three terms in Equations \hyperref[Eq_TTS1]{\ref*{Eq_TTS1}}, \hyperref[Eq_TAC]{\ref*{Eq_TAC}} and \hyperref[Eq_TC]{\ref*{Eq_TC}}.
Observe that the terms $TTT$ and $TAC$ are linear, that is,

\begin{equation}
TTT = \int_0^{t_f} c \cdot N(t) \mathrm{dt}    
\end{equation}

\begin{equation}
TAC = \sum_{t_a \in \mathcal{T}_a} \int_0^{t_f} s \cdot Q^{t_a}(N(t),\gamma(t)) \cdot \mathcal{L}(t_a,t) \mathrm{dt}     
\end{equation}

\noindent We take the terminal criterion $TC$ as quadratic:

\begin{equation}
TC = \frac{1}{2}\cdot N(t_f) \cdot \mu \cdot N(t_f)    
\end{equation}

The vector (or matrix) $\mu$ in Equation~\eqref{Eq_TC} represents the weights assigned to the final accumulation of vehicles in each region of the network. It is used to penalize any remaining vehicles at the end of the time horizon $t_f$, thereby encouraging the system to clear the network efficiently. If $\mu$ is taken as a diagonal matrix with entries $\mu_i$, the terminal cost becomes:

\begin{equation}
TC = \frac{1}{2} \sum_{i \in \mathcal{J} \cup \mathcal{O}} \mu_i \cdot N_i^2(t_f),
\end{equation}

\noindent which penalizes each region proportionally to the square of its residual accumulation. Larger values of $\mu_i$ imply a stronger preference for minimizing the number of vehicles in region $i$ at the terminal time. If $\mu=0$, the final state is unpenalized. Thus, the choice of $\mu$ influences how aggressively the system seeks to empty the network before $t_f$.

\noindent All in all, the problem is formulated as an optimal control problem:

\begin{equation}\label{eq_DefCriterion}
\min \int_0^{t_f} \Bigg[c \cdot N(t) + \sum_{t_a \in \mathcal{T}_a} s \cdot Q^{t_a}(N(t),\gamma(t)) \cdot \mathcal{L}(t_a,t)\Bigg] \mathrm{dt} + \frac{1}{2}\cdot N(t_f) \cdot \mu \cdot N(t_f)    
\end{equation}

\noindent subject to constraints defined in the set of equations \hyperref[Eq_Combined]{\ref*{Eq_Combined}}.

\subsection{Calculation of the gradient}\label{Sect3.3}
This section discusses how we calculate the gradient to solve the previous optimal control formulation. Considering that the initial state is given, the criterion can be viewed as a function of the decision variables only. We consider the time steps $k=0,\cdots,K$. Furthermore, we simplify some notation to derive the gradient of the criterion concerning the decision variables in a discrete-time setting. 

Let us start by introducing some notation. We define: $N_k$ as the state in the time step $k$; $\beta_k$ as the command at the time step $k$ ($\beta_k = (\gamma_k,\Delta_k)$); $g_k(N_k,\beta_k)$ as the contribution of the time step $k$ to the criterion; and $g_K(N_K)$ as the terminal criterion. Therefore:

\begin{equation}
g_k(N_k,\beta_k) = \Bigg[c \cdot N_k + \sum_{t_a \in \mathcal{T}_a} s \cdot Q^{t_a}(N_k,\beta_k) \cdot \mathcal{L}(t_a,k\cdot \Delta t) \Bigg] \times \Delta t    \quad\forall k=0..K-1    
\end{equation}

\noindent where vector $c$ represents the weights applied to the number of vehicles present in each region of the network, thus capturing the Total Time Spent (TTS) component in the criterion. Each component $c_i$ reflects the marginal cost (or disutility) of vehicle accumulation in region $i$. If $c$ is taken as a vector of ones, the TTS corresponds to the total number of vehicles in the network at each time step. Alternatively, $c$ can be region-dependent to prioritize congestion mitigation in specific areas. Vector $s$ represents the weights assigned to the early or late arrival penalty for each flow component in $Q^{t_a}$. It scales the penalty $\mathcal{L}(t_a,t)$ incurred by travellers who arrive at time $t$ with a desired arrival time $t_a$. Each component of $s$ reflects the importance or sensitivity of arrival punctuality, and can be used to prioritize critical routes or destinations. The term $s \cdot Q^{t_a}(N_k,\beta_k) \cdot \mathcal{L}(t_a,k\cdot \Delta t)$ thus captures the total time-dependent disutility due to arrival mismatches across the network.

\begin{equation}
g_K(N_K) = \frac{1}{2}\cdot N_K \cdot \mu \cdot N_K    
\end{equation}

\begin{equation}\label{Eq_State_1}
N_{k+1} = E(N_k,\beta_k)     \quad\forall k=0..K-1    
\end{equation}

\begin{equation}\label{Eq_State_2}
N_{k+1} = N_k + \Delta t \cdot F(N_k,\beta_k) 
\end{equation}

\noindent i.e.:

\begin{equation}\label{Eq_State_3}
E(N_k,\beta_k) = \Delta t \cdot F(N_k,\beta_K) + N_k \quad\forall k=0..K-1    
\end{equation}

\noindent Equations \hyperref[Eq_State_1]{\ref*{Eq_State_1}}, \hyperref[Eq_State_2]{\ref*{Eq_State_2}}, \hyperref[Eq_State_3]{\ref*{Eq_State_3}} constitute the discretised state equations of the system. We define $J$ as the discretised version of criterion \hyperref[eq_DefCriterion]{\ref*{eq_DefCriterion}}:

\begin{equation}
J \stackrel{\text{def}}{=} \sum_{k=0}^{K-1} g_k(N_k,\beta_k) + g_K(N_K)    
\end{equation}

The generic calculation of the derivatives of $J$ concerning $\beta$, i.e. $\partial_{\beta_{K-1}}J$ and $\partial_{\beta_{k}}J$. First, note that $N_h$ depends on $\beta_k$ if and only if $h>k$. More precisely the following apply
\begin{align}
\begin{array}{lcll}
     \partial_{\beta_k} N_{k+1} & = & \partial_{\beta_k} E_{k} & \\
     \partial_{\beta_k} N_{h+1} & = & \partial_{N_h} E_{h} \cdot \partial_{\beta_k} N_{h} & \forall h>k
\end{array}
\end{align}
Second, for all $k$:
\begin{align}
\partial_{\beta_{k}}J = \partial_{\beta_{k}}g_{k} + \sum_{h=k+1}^{K-1} \partial_{N_h}g_h \cdot \partial_{\beta_{k}}N_h + \partial_{N_K}g_K \cdot \partial_{\beta_{h}}N_K 
\end{align}
\noindent and the iterative calculation of the derivatives of $J$ concerning $\beta_k$ follows:
\begin{equation}
\begin{array}{lcl}
\partial_{\beta_{K-1}}J & = & \partial_{\beta_{K-1}}g_{K-1}  + \partial_{N_K}g_K \cdot \partial_{\beta_{K-1}}E_{K-1} \\
\partial_{\beta_{k}}J & = & \partial_{\beta_{k}}g_{k} + \sum_{h=k+1}^{K} \partial_{N_h}g_h \cdot \partial_{N_{h-1}} E_{h-1} \cdot \ldots \cdot \partial_{N_{k+1}} E_{k+1} \cdot \partial_{\beta_{k}}E_{k} \quad \forall k<K
\end{array}
\end{equation}

Thus, let us introduce the adjoint state $\Upsilon$:
\begin{equation}\label{Eq_Upsilons}
\begin{cases}
\Upsilon_K = \partial_{N_{K}}g_K \\
\Upsilon_k = \partial_{N_{k}}g_k + \Upsilon_{k+1} \cdot \partial_{N_{k}} E_k, \quad k=K-1, \cdots, 0
\end{cases}
\end{equation}
and the gradient of $J$ is expressed as:
\begin{equation}\label{Eq_partialbeta}
\partial_{\beta_{k}}J = \partial_{\beta_{k}}g_k + \Upsilon_{k+1} \cdot \partial_{\beta_{k}} E_k , \quad k=0,\dots,K
\end{equation}

Actually, given that the set of constraints applied to $\beta$ is convex (constraints defined in Equations \hyperref[Eq_TotalDemand1]{\ref*{Eq_TotalDemand1}}, \hyperref[Eq_Constraint1]{\ref*{Eq_Constraint1}}, \hyperref[Eq_Constraint2]{\ref*{Eq_Constraint2}}, \hyperref[Eq_Constraint3]{\ref*{Eq_Constraint3}} and \hyperref[Eq_Constraint4]{\ref*{Eq_Constraint4}}), the discrete-time minimum principle applies. We will not apply this approach here.

\subsection{Optimization scheme (generic aspects)}\label{Sect3.4}
We propose an iterative optimisation scheme based on the projected gradient as follows:

\begin{equation}\label{Eq_ConjugateGradient}
\beta^{T+1} = \mathcal{P}_{\mathcal{K}}\left[ \beta^\tau  - \alpha^\tau  \cdot \partial_{\beta}J(\beta^\tau ) \right]  
\end{equation}

\noindent where $\tau $ denotes the iteration index, $\mathcal{P}_{\mathcal{K}}$ is the projector on $\mathcal{K}$, the set of constraints applying to $\beta$; and $\alpha^\tau $ satisfies the divergent series rule of $\alpha^\tau \rightarrow 0$ as $\tau  \rightarrow \infty$, and $\sum_\tau  \alpha^\tau = \infty$.

For calculating $\partial_{\beta}J(\beta^\tau )$, we apply the equations \hyperref[Eq_Upsilons]{\ref*{Eq_Upsilons}} and \hyperref[Eq_partialbeta]{\ref*{Eq_partialbeta}}, and the state equation:

\begin{equation}
N_{k+1}^\tau = E(N_k^\tau,\beta_k^\tau), k=0,\dots,K-1    
\end{equation}

Considering that the calculation of $\partial_{\beta}J(\beta)$ is quite complicated, and the convergence of Equation \hyperref[Eq_ConjugateGradient]{\ref*{Eq_ConjugateGradient}} is slow, we propose to start with a time step large enough. Let Equation \hyperref[Eq_ConjugateGradient]{\ref*{Eq_ConjugateGradient}} diminish the criterion, then divide the time step into two and restart iteration \hyperref[Eq_ConjugateGradient]{\ref*{Eq_ConjugateGradient}} with as initial $\beta$ the command obtained with the larger step solution.

For the initial assignment, we could:

\begin{itemize}
    \item consider for each region $i$, an average travel time $\theta_i = 1 / \partial_{N_i}\Delta_i(N_i)_{|N_i=0}$;
    \item assign traffic to the regional paths with the minimal travel times connecting the Origin-Destination pair $\forall w = (O,D) \in \mathcal{W}$, yielding the constant coefficients $\gamma_{ij}^{D,t_a}$. The path travel times are calculated based on the $\theta_i $;
    \item assign the departure times per Origin-Destination pair in such a way that traffic arrives at the Destination as close as possible to $t_a$, yielding the demand $\Delta_O^{D,t_a}(t)$ for $t=k\cdot \Delta t$.
\end{itemize}

The question at hand is how to establish the suitable time interval, denoted $\Delta t$. We suggest utilizing:

\begin{equation}
\begin{split}
\Delta t_{init.} &= \min_{i \in \mathcal{J}} \theta_i \\
 &= \min_{i \in \mathcal{J}} 1 / \partial_{N_i} \Delta_i(N_i)_{|N_i=0}
\end{split}    
\end{equation}

\noindent $\theta_i$  can be interpreted as the minimum travel time of region $i$, refer to \citep{khoshyaran2025dynamic}, subsubsection 2.2.5. Considering the choice of $\Delta_O(.), \forall O \in \mathcal{O}$, we should pay special attention to the dynamics of $N_O, \forall O \in \mathcal{O}$.

To ensure numerical stability in the discrete-time simulation of the system-optimal control dynamics, a Courant–Friedrichs–Lewy (CFL) condition is required. This constraint guarantees that no traveller may traverse more than one region per time step, thus preserving the physical consistency of the regional propagation model.










\subsubsection{Numerical Stability and CFL Conditions}\label{Sect3.cfl}

The discretisation necessitates a CFL-type condition, ensuring that within the discretised model, no traveller can traverse more than one region per time step. Moreover, the time step should be maximised to minimise numerical viscosity.

The following expression results for the time step $\Delta t$ (refer again to \citep{khoshyaran2025dynamic}):

\begin{equation}
\Delta t = \min \Bigg(\min_{i \in \mathcal{I}\cup\mathcal{O}}(\partial_{N_i}\Delta_i(N_i))^{-1},\min_{j \in \mathcal{I}}(-\partial_{N_j}\Sigma_j(N_j))^{-1}\Bigg)    
\end{equation}

In particular, we must take into consideration the origin cells (i.e. the special case of origins described in Section \hyperref[Sect2.2.1]{\ref*{Sect2.2.1}}). In particular, $\mu_O$ must be chosen in such a way that $\partial_{N_O} \Delta_O = \frac{\mu_O}{Q_{Ox}}$ is not too small. Otherwise, we would need to choose a value of $\Delta t$ which would be small, resulting in numerical smoothing of all dynamical variables. On the other hand, if $\mu_O$ is chosen large, there is a significant smoothing effect on the demands $\Delta_O^{D,t_a}(t)$. Thus, a compromise is required.

\section{Optimization of J}\label{Sect4}
We apply Equations \hyperref[Eq_ConjugateGradient]{\ref*{Eq_ConjugateGradient}}, \hyperref[Eq_Upsilons]{\ref*{Eq_Upsilons}} and \hyperref[Eq_partialbeta]{\ref*{Eq_partialbeta}}, therefore we need to calculate the gradients $\partial_{\beta_{k}}g_{k}$, $\partial_{\beta_{k}}E_{k}$, $\partial_{N_{k}}E_{k}$. In the sequel, we only calculate the non-zero terms. Let us describe the coordinates of the variables:

\begin{itemize}
    \item For $N_k$ (state variables), we consider $N_i(k), \forall i \in \mathcal{O} \cup \mathcal{J}$ and $N_i^{D,t_a}(k), \forall i \in \mathcal{O} \cup \mathcal{J} \wedge \forall D \in \mathcal{D} \wedge t_a \in \mathcal{T}_a$. Note that we write $k$ where we should write $k\cdot \Delta t$.
    \item For $\Upsilon_k$ (the adjoint state), we consider $\Upsilon_i(k)$ and $\Upsilon_i^{D,t_a}(k)$ with the same range as for $N_k$.
    \item For $\beta_k$, we consider $\gamma_{ij}^{D,t_a}(k), \forall i \in \mathcal{O} \cup \mathcal{J} \wedge j \in \Gamma^+(i) \wedge D \in \mathcal{D} \wedge t_a \in \mathcal{T}_a$ and $\Delta^{D,t_a}_O(k), \forall O \in \mathcal{O} \wedge D \in \mathcal{D} \wedge t_a \in \mathcal{T}_a$. 
\end{itemize}

The solution procedure consists of four steps.

\vspace{2ex}
\noindent \textbf{Step 1:}

The \textbf{first step} of the optimisation procedure consists of the initialisation of $\Upsilon$ following equation \hyperref[Eq_Upsilons]{\ref*{Eq_Upsilons}}:
 
\begin{equation}
\left| \; \begin{array}{l} 
   \Upsilon_i(K) = \mu_i \cdot N_i(K)      \\
   \Upsilon_i^{D,t_a}(K) = 0.     
\end{array} \right.
\end{equation}

\vspace{2ex}
\noindent \textbf{Step 2:}

The \textbf{second step} consists of calculating the derivatives of the $g_k$ and $E_k$ terms related to the state variables $N_k$ to evaluate the adjoint variables defined in \hyperref[Eq_Upsilons]{\ref*{Eq_Upsilons}}. For this, we need to determine the partial derivatives of $g_k$ with respect to $N_k$ and $\beta_k$. The difficulty lies in the fact that we need to chain derivatives, since the calculations of relevant quantities are chained. 

We start by considering the outflow of the network at a destination $D$. In Equation \hyperref[Eq_Supplyij]{\ref*{Eq_Supplyij}} we consider $i\in\Gamma^{-}(D)$ and set $\ell = D$ in order to express that traffic bound to $D$ exits in $D$. Hence for $i\in\Gamma^{-}(D)$, $\gamma_{iD}^{D,t_a}=1$, $\gamma_{iD}^{d,t_a}=0$ if $d\not= D$ and Equation \hyperref[Eq_Supplyij]{\ref*{Eq_Supplyij}} reduces to:

\begin{equation}\label{Eq_deltaiD}
\delta_{iD}(t) = \delta_i(t) \cdot \frac{\sum_{t_a \in \mathcal{T}_a} N_i^{D,t_a}(t)}{N_i(t)}, \quad \forall D \in \mathcal{D} \wedge  \forall i \in \Gamma^-(D) 
\end{equation}
Equation \hyperref[Eq_Inflows1]{\ref*{Eq_Inflows1}} yields:  
\begin{equation}\label{Eq_qiD}
q_{iD}(t) = \min\left[\sigma_{iD}(t),\delta_{iD}(t)\right]    
\end{equation}
Further, $q_{iD}^{D,t_a}(t)$ is proportional to the composition of traffic concerning the desired arrival time, thus:
\begin{equation}\label{Eq_qiDt_a}
q_{iD}^{D,t_a}(t) =  q_{iD}(t) \frac{N_i^{D,t_a}(t)}{\sum_{\theta \in \mathcal{T}_a} N_i^{D,\theta}(t)}  
\end{equation}

We deduce expressions for the flows $Q^{D,t_a}_{iD}$. The demand $\sigma_{iD}(k)$ is given from the data.  We introduce a dummy variable $s_{iD}$ which equals 1 if $q_{iD}=\sigma_{iD}$, and 0 if $q_{iD}=\delta_{iD}$ (i.e. $q_{iD}<\sigma_{iD}$). Note that this is not a new variable. The main purpose is only to simplify the calculations. 

Changing $\delta_i(t)$ with $\Delta_i (N_i)$ and substituting Equation \hyperref[Eq_deltaiD]{\ref*{Eq_deltaiD}} into Equation \hyperref[Eq_qiD]{\ref*{Eq_qiD}} and Equation \hyperref[Eq_qiD]{\ref*{Eq_qiD}} into Equation \hyperref[Eq_qiDt_a]{\ref*{Eq_qiDt_a}}, we obtain the following.

\begin{align}\label{Eq_qiD_with_sid}
\begin{split}
q_{iD}^{D,t_a} &= \min\Bigg[\sigma_{iD} \cdot \frac{N_i^{D,t_a}}{\sum_{\theta \in \mathcal{T}_a} N_i^{D,\theta}}, \Delta_i(N_i) \cdot \frac{N_i^{D,t_a}}{N_i}\Bigg] \\
 &= s_{iD} \cdot \sigma_{iD} \cdot \frac{N_i^{D,t_a}}{\sum_{\theta \in \mathcal{T}_a} N_i^{D,\theta}} + (1-s_{iD}) \cdot \Delta_i(N_i) \cdot \frac{N_i^{D,t_a}}{N_i}
\end{split}
\end{align}

We calculate the derivative terms of $q_{iD}^{D,t_a}$ concerning $N_i$, $N_i^{D,t_a}$ and $N_i^{D,\theta}$. We obtain the following.

\begin{equation}\label{Eq_PartialqiD1}
\partial_{N_i} q_{iD}^{D,t_a} = \begin{cases}
0 & \textrm{if } s_{iD} = 1 \\
\Delta_i^{'}(N_i) \cdot \frac{N_i^{D,t_a}}{N_i} - \Delta_i(N_i) \cdot \frac{N_i^{D,t_a}}{N_i^2} & \textrm{if } s_{iD} = 0
\end{cases}
\end{equation}

\begin{equation}\label{Eq_PartialqiD2}
\partial_{N_i^{D,t_a}} q_{iD}^{D,t_a} = \begin{cases} \frac{\sigma_{iD}}{\sum_{\theta \in \mathcal{T}_a} N_i^{D,\theta}} - \frac{\sigma_{iD} \cdot N_i^{D,t_a}}{\left(\sum_{\theta \in \mathcal{T}_a} N_i^{D,\theta}\right)^2} & \textrm{if } s_{iD} = 1 \\
\frac{\Delta_i(N_i) }{N_i} & \textrm{if } s_{iD} = 0
\end{cases}
\end{equation}

\begin{equation}\label{Eq_PartialqiD3}
\partial_{N_i^{D,\theta}} q_{iD}^{D,t_a} = \begin{cases}
- \frac{\sigma_{iD} \cdot N_i^{D,t_a}}{\left(\sum_{\theta \in \mathcal{T}_a} N_i^{D,\theta}\right)^2} & \textrm{if } s_{iD} = 1 \\
0 & \textrm{if } s_{iD} = 0
\end{cases}, \textrm{for } \theta \neq t_a
\end{equation}

Finalisation of the calculation of $\partial_{N_k}g_k$:

\begin{equation}\label{Eq_Partialgk1}
\partial_{N_i(k)} g_k = \Biggl(1 + \sum\limits_{\substack{D \in \Gamma^+(i) \cap \mathcal{D} \\ t_a \in \mathcal{T}_a}} \partial_{N_i} q_{iD}^{D,t_a}(k) \cdot \mathcal{L}(t_a,k\cdot \Delta t)\Biggl)\cdot \Delta t
\end{equation}

\begin{equation}\label{Eq_Partialgk2}
\partial_{N_i^{Dd,t_a}(k)} g_k = \Biggl(\sum\limits_{\substack{D \in \Gamma^+(i) \cap \mathcal{D} \\ t_a \in \mathcal{T}_a}} \partial_{N_i^{D,t_a}}q_{iD}^{D,\theta} \cdot \mathcal{L}(\theta,k\cdot \Delta t)  \Biggl) \cdot \Delta t   
\end{equation}

These two equations summarise the calculation of the partial derivative $\partial_{N_k}g_k$. Note that Equations \hyperref[Eq_PartialqiD1]{\ref*{Eq_PartialqiD1}}, \hyperref[Eq_PartialqiD2]{\ref*{Eq_PartialqiD2}}, \hyperref[Eq_PartialqiD3]{\ref*{Eq_PartialqiD3}}, \hyperref[Eq_Partialgk1]{\ref*{Eq_Partialgk1}}, and \hyperref[Eq_Partialgk2]{\ref*{Eq_Partialgk2}} also show that $\partial_{N_k}g_k$ is $O$ unless there is a couple $(OD) \in \mathcal{W}$ such that $D \in \Gamma^+(O)$.

\vspace{2ex}
\noindent \textbf{Step 3:}

In the \textbf{third step}, we calculate the partial derivative of $E_k$ with respect to $N_k$. For this, we first solve Equations \hyperref[Eq_DNiDta]{\ref*{Eq_DNiDta}} and \hyperref[Eq_DNi]{\ref*{Eq_DNi}} using the time discretisation as follows:

\begin{equation}\label{Eq_Discrek1}
N_i^{D,t_a}(k+1) = N_i^{D,t_a}(k) + \Bigg[\sum_{l \in \Gamma^-(i)} q_{li}^{D,t_a}(k)-\sum_{j \in \Gamma^+(i)} q_{ij}^{D,t_a}(k) \Bigg]\cdot \Delta t := E_i^{D,t_a}(k)
\end{equation}

\begin{equation}\label{Eq_Discrek2}
N_i(k+1) = N_i(k) + \Bigg[\sum_{l \in \Gamma^-(i)} q_{li}(k)-\sum_{j \in \Gamma^+(i)} q_{ij}(k)\Bigg]\cdot \Delta t := E_i(k)
\end{equation}

The right-hand side of Equations \hyperref[Eq_Discrek1]{\ref*{Eq_Discrek1}} and \hyperref[Eq_Discrek2]{\ref*{Eq_Discrek2}} constitutes $E_k$. Again, we must express chainwise $E_k$ as a function of $\beta_k$ and $N_k$. For this, we start with Equation \hyperref[Eq_qijDta]{\ref*{Eq_qijDta}}:

\begin{equation}
q_{ij}^{D,t_a}(k) = q_{ij}(k) \cdot \gamma_{ij}^{D,t_a} \cdot \frac{N_i^{D,t_a}(k)}{\sum_{\substack{\delta \in \mathcal{D} \\ \theta \in \mathcal{T}_a}} \gamma_{ij}^{\delta,\theta}(k) \cdot N_{i}^{\delta,\theta}(k)}   
\end{equation}

\noindent where (STRADA model \hyperref[Eq_DemandSTRADA]{\ref*{Eq_DemandSTRADA}}):

\begin{equation}
q_{ij}(k) = \min \left[ \beta_{ij}\cdot \Sigma_j(N_j(k)), \delta_{ij}(k) \right]   
\end{equation}

\noindent and,

\begin{equation}
\delta_{ij}(k) = \Delta_i(N_i(k)) \cdot \sum_{\theta,\delta} \gamma_{ij}^{\theta,\delta} \cdot \frac{N_i^{\delta,\theta}(k)}{N_{i}(k)}  
\end{equation}

Note that if $i \in \mathcal{O}$, we apply Equations \hyperref[Eq_DNODta]{\ref*{Eq_DNODta}} and \hyperref[Eq_DNO]{\ref*{Eq_DNO}}, that is: 

\begin{equation}\label{Eq_NODtak1}
N_{O}^{D,t_a}(k+1) = N_O^{D,t_a}(k) + \Delta t \cdot \Bigg[\Delta_O^{D,t_a}(k) - \sum_{i \in \Gamma^+(O)} q_{Oi}^{D,t_a}(k)\Bigg] := E_O^{D,t_a}   
\end{equation}

\begin{equation}\label{Eq_NODtak}
N_{O}(k+1) = N_{O}(k) + \Delta t \cdot \Bigg[\sum_{\delta,\theta} \Delta_O^{\delta,\theta}(k) - \sum_{i \in \Gamma^+(O)} q_{Oi}(k)\Bigg] := E_O
\end{equation}

These two equations are similar to Equations \hyperref[Eq_Discrek1]{\ref*{Eq_Discrek1}} and \hyperref[Eq_Discrek2]{\ref*{Eq_Discrek2}}, but exhibit a different structure concerning the $\beta$ dependency, i.e. they contains the command variables $\Delta_O^{D,t_a}$. 

Calculating the derivatives of Equation \hyperref[Eq_Discrek2]{\ref*{Eq_Discrek2}} concerning $N_i$:

\begin{equation}\label{Eq_partialEikNik}
\partial_{N_i(k)} E_i(k) = 1 + \Delta t \cdot \Bigg[\sum_{\ell \in \Gamma^-(i)} \partial_{N_i (k)} q_{\ell i}(k) - \sum_{j \in \Gamma^+(i) }  \partial_{N_i (k)} q_{ij}(k)  \Bigg]   
\end{equation}
 
\begin{equation}\label{Eq_partialEikNlk}
\partial_{N_\ell (k)} E_i(k) = \Delta t \cdot \partial_{N_\ell (k)} q_{\ell i}(k), \quad \forall \ell  \in \Gamma^-(i) 
\end{equation}

\noindent and,

\begin{equation}\label{Eq_partialEikNjk}
\partial_{N_j(k)} E_i(k) = -\Delta t \cdot \partial_{N_j(k)} q_{ij}(k), \quad\forall j \in \Gamma^+(i)   
\end{equation}

Here, we have a specific point to attend to, since a state depends on neighbouring states for its dynamics.

We recall the following rule: All derivatives not explicitly calculated are null. 

We also calculate other derivatives. For this, consider $u,v \in \mathcal{J}$. The inflow from the region $u$ to $v$ is $q_{uv}=\min\left[\beta_{uv}\cdot \Sigma_v(N_v),\delta_{uv}\right]$. We also define $s_{uv}$ as a binary variable equal to 1 if $q_{uv}=\beta_{uv}\cdot \Sigma_v(N_v)$ and 0 if $q_{uv}=\delta_{uv} < \beta_{uv}\cdot \Sigma_v(N_v)$. This is a standing notation that is required for the calculation of the derivatives. 

\begin{equation}
\partial_{N_i (k)} q_{li} (k) = \begin{cases}
\beta_{\ell i} \Sigma_i^{'}(N_i (k)) & \textrm{if } s_{\ell i} = 1\\
\partial_{N_i (k)} \delta_{\ell i (k)} = 0 & \textrm{if } s_{\ell i} = 0
\end{cases} 
\end{equation}

\noindent where $\delta_{\ell i} (k)$ depends on the states $N_\ell  (k)$.

\begin{equation}
\partial_{N_i} (k) q_{ij} (k) = \begin{cases}
0 & \textrm{if } s_{li} = 1\\
\partial_{N_i (k)} \delta_{ij} (k) & \textrm{if } s_{ij} = 0
\end{cases}    
\end{equation}

\noindent where $\sigma_{ij} (k)$ depends only on $N_j (k)$. 

Recall, generically:

\begin{equation}
\delta_{uv} = \Delta_u(N_u) \cdot \sum_{\theta,\delta} \gamma^{\theta,\delta}_{uv} \cdot \frac{N_u^{\theta,\delta}}{N_u}
\end{equation}

The partial derivative of $\delta_{uv}$ regarding $N_u$ is:

\begin{align}
\begin{split}
\partial_{N_u} \delta_{uv} &= \Delta^{'}(N_u) \cdot \sum_{\theta,\delta} \gamma^{\theta,\delta}_{uv} \cdot \frac{N_u^{\theta,\delta}}{N_u} - \Delta_u(N_u) \cdot \sum_{\theta,\delta} \gamma^{\theta,\delta}_{uv} \cdot \frac{N_u^{\theta,\delta}}{N_u^2} \\
&= \Bigg(\frac{\Delta^{'}_u(N_u)}{\Delta_u(N_u)}-\frac{1}{N_u} \Bigg)\cdot \delta_{uv}
\end{split}
\end{align}

Therefore:

\begin{equation}
\partial_{N_i (k)} \delta_{ij} (k) = \Bigg(\frac{\Delta^{'}_i(N_i)}{\Delta_i(N_i)}-\frac{1}{N_i} \Bigg)\cdot \delta_{ij} \quad \left( \; = \partial_{N_i (k)} q_{ij} (k)  \mbox{ if   } s_{ij} = 0 \, \right)
\end{equation}

This completes the calculation of Equation \hyperref[Eq_partialEikNik]{\ref*{Eq_partialEikNik}}.

Also, for $\ell \in \Gamma^-(i)$, from Equation \hyperref[Eq_partialEikNlk]{\ref*{Eq_partialEikNlk}}, the partial derivative of $E_i (k)$ concerning $N_\ell (k)$ is:

\begin{align}
\begin{split}
\partial_{N_\ell (k)} E_i (k) &= \Delta t \cdot \partial_{N_\ell (k)} q_{\ell i} (k) \\
 &= \begin{cases}
 0 & \textrm{if } s_{\ell i} = 1 \\
 \Delta t \cdot \Bigg(\frac{\Delta_\ell^{'} (N_\ell (k))}{\Delta_\ell (N_\ell (k))}-\frac{1}{N_\ell (k)} \Bigg)\cdot \delta_{\ell i} (k) & \textrm{if } s_{\ell i} = 0
 \end{cases}
\end{split}
\end{align}

Similarly, for $j \in \Gamma^+(i)$, from Equation \hyperref[Eq_partialEikNjk]{\ref*{Eq_partialEikNjk}}, the partial derivative of $E_i$ regarding $N_j$ is:

\begin{align}
\begin{split}
\partial_{N_j (k)} E_i (k) &= -\Delta t \cdot \partial_{N_j (k)} q_{ij} (k) \\
&= \begin{cases}
 -\Delta t \cdot \beta_{ij} \cdot \Sigma^{'}_j(N_j (k)) & \textrm{if } s_{ij} = 1 \\
 0 & \textrm{if } s_{li} = 0
 \end{cases}
\end{split}
\end{align}


The Origin regions constitute special cases, i.e. when $i = O \in \mathcal{O}$. All the previous calculations for the partial derivatives are similar for the Origin regions. There is only one main difference: $\Gamma^- (O)$ is empty. Therefore, there are no summation terms such as $\sum_{j \in \Gamma^- (O)} ...$, and there are no derivatives concerning $N_\ell$'s, $\ell \in \Gamma^- (O)$.

\vspace{2ex}
\noindent \textbf{Step 4:}

In the \textbf{fourth step}, we calculate the derivatives concerning $\beta_k$. First, we emphasize that $\partial_{\beta_k} g_k = 0$, since the only terms which could have a non-zero contribution are the $Q^{t_a}$, but these do not depend on the $\gamma$'s since $\gamma_{iD}^{D,t_a}=1, \forall D, i \in \Gamma^{-}(D) \wedge t_a \in \mathcal{T}_a$. 

The calculation of the partial derivative of $E_k$ concerning $\beta_k$, i.e. $\partial_{\beta_k} E_k$. First, we focus on the Origin regions. 

Recall that:

\begin{align}
&N_{O}^{D,t_a}(k+1) = N_O^{D,t_a}(k) + \Delta t \cdot \Bigg[\Delta_O^{D,t_a}(k) - \sum_{i \in \Gamma^+(O)} q_{Oi}^{D,t_a}(k)\Bigg] \nonumber \\
&N_{O}(k+1) = N_{O}(k) + \Delta t \cdot \Bigg[\Delta_O(k) - \sum_{i \in \Gamma^+(O)} q_{Oi}(k)\Bigg] \nonumber  \\
&q_{Oi}^{D,t_a}(k) = q_{Oi}(k) \cdot \gamma_{Oi}^{D,t_a}(k) \cdot \frac{N_O^{D,t_a}(k)}{\sum_{\substack{\theta \in \mathcal{D} \\ \theta \in \mathcal{T}_a}} \gamma_{Oi}^{\delta,\theta}\cdot N_O^{\delta,\theta}(k)} \\
&q_{Oi}(k) = \min(\beta_{Oi} \cdot \Sigma_i(N_i(k)), \delta_{Oi}(k)) \nonumber  \\
&\delta_{Oi}(k) = \Delta_O(N_O(k)) \cdot \sum_{\substack{\theta \in \mathcal{D} \\ \theta \in \mathcal{T}_a}} \gamma_{Oi}^{\delta,\theta}(k) \cdot \frac{N_0^{\delta,\theta}(k)}{N_O(k)} \nonumber 
\end{align}






We can calculate here the only nonzero derivatives concerning $\Delta_O$ and $\Delta_O^{D,t_a}$:

\begin{equation}
 \partial_{\Delta_O^{D,t_a}(k)} E_O^{D,t_a} (k) = \partial_{\Delta_O (k)} E_O (k) = \Delta t  
\end{equation}

Recall the constraints $\Delta_O(k) - \sum_{\delta,\theta} \Delta_O^{\delta,\theta} (k) = 0, \; \forall O \in \mathcal{O} \wedge \forall k = 1,\dots,K-1$. Based on this constraint, we eliminate $\Delta_O(k)$ and replace $\Delta_O(k)$ with $\sum_{\theta,\delta} \Delta_O^{\delta,\theta}(k)$ hereafter. 

The next step consists of calculating the gradient of $E_i^{D,t_a}$. We can also ascertain $\partial_{\gamma_{h\ell}^{\delta,\theta}} E_i^{D,t_a}=0, \mathrm{if } h\neq i \wedge \ell \neq i$. Let us recall the following equations:

\begin{equation}
E_i^{D,t_a}(k) = N_i^{D,t_a}(k) + \Bigg[\sum_{l \in \Gamma^-(i)} q_{li}^{D,t_a}(k)-\sum_{j \in \Gamma^+(i)} q_{ij}^{D,t_a}(k) \Bigg]\cdot \Delta t    
\end{equation}

\begin{equation}
q_{uv}(k)=\min[\beta_{uv}\cdot \Sigma_v(N_v(k)),\delta_{uv}(k)]    
\end{equation}

\begin{equation}
\delta_{uv}(k) = \Delta_u(N_u(k)) \cdot \sum_{\theta,\delta} \gamma^{\theta,\delta}_{uv}(k) \cdot \frac{N_u^{\theta,\delta}(k)}{N_u(k)}
\end{equation}

\noindent and,

\begin{equation}
q_{uv}^{D,t_a}(k) = q_{uv}(k) \cdot \gamma_{uv}^{D,t_a}(k) \cdot \frac{N_u^{D,t_a}(k)}{\sum_{\delta,\theta} \gamma_{uv}^{\delta,\theta}(k) \cdot N_u^{\delta,\theta}(k)}    
\end{equation}

For immediate use, we reformulate the expression of $q_{uv}^{D,t_a}$:

\begin{equation}\label{Eq_quvDta}
q_{uv}^{D,t_a}(k) =\min\Bigg[\beta_{uv}\cdot \Sigma_v(N_v(k)) \cdot \frac{\gamma_{uv}^{D,t_a}(k) \cdot N_u^{D,t_a}(k)} {\sum_{\delta,\theta} \gamma_{uv}^{\delta,\theta}(k) \cdot N_u^{\delta,\theta}(k)}, \Delta_u(N_u(k))\cdot \frac{\gamma_{uv}^{D,t_a}(k)\cdot N_u^{D,t_a}(k)}{N_u(k)} \Bigg]    
\end{equation}

\noindent As before, we define $s_{uv}$ as follows. $s_{uv}=1$ if $q_{uv}^{D,t_a}$ equals supply, and $s_{uv}=0$ if $q_{uv}^{D,t_a}$ equals demand within the minimum operator of \hyperref[Eq_quvDta]{\ref*{Eq_quvDta}}.

Given the Equation \hyperref[Eq_quvDta]{\ref*{Eq_quvDta}}, it follows that the only non-zero derivatives of $q_{uv}^{D,t_a}$ concerning $\gamma$ are obtained concerning $\gamma_{uv}^{D,t_a}$. Therefore, we deduce the following expression for the partial derivatives of the states $E_i^{D,t_a}(k)$ concerning $\gamma(k)$:

\begin{equation}\label{Eq_DerEwrtGamma}
\begin{cases}
\partial_{\gamma_{hi}^{\delta,\theta}(k)} E_i^{D,t_a}(k) = \Delta t \cdot \partial_{\gamma_{hi}^{\delta,\theta}(k)} q_{hi}^{D,t_a}(k) \\
\partial_{\gamma_{i\ell}^{\delta,\theta}(k)} E_i^{D,t_a}(k) = -\Delta t \cdot \partial_{\gamma_{i\ell}^{\delta,\theta}(k)} q_{i\ell}^{D,t_a}(k)
\end{cases}
\end{equation}

The partial derivatives of $q_{hi}^{D,t_a}(k)$ concerning $\gamma_{hi}^{\delta,\theta}(k)$ is deduced from (\ref{Eq_quvDta}):
\begin{equation}
    \partial_{\gamma_{hi}^{\delta,\theta}(k)} q_{hi}^{D,t_a}(k) = \; \left\{ \;  
        \begin{array}{lcl}
            \beta_{hi} \cdot \Sigma_i(N_i (k)) \cdot \frac{ N_h^{D,t_a}(k)}{\sum_{\delta',\theta'} \gamma_{hi}^{\delta',\theta'}(k) \cdot N_h^{\delta',\theta'}(k)} \cdot\left[ 1 - \frac{N_h^{D,t_a}(k)\cdot \gamma_{hi}^{D,t_a}(k)}{\sum_{\delta',\theta'} \gamma_{hi}^{\delta',\theta'}(k) \cdot N_h^{\delta',\theta'}(k)} \right]  & \textrm{if } &  s_{hi}=1 \wedge \delta, \theta = D, t_a \\
            -\beta_{hi}\cdot \Sigma_i(N_i (k))\cdot \frac{N_h^{D,t_a}(k) \cdot N_h^{\delta,\theta}(k) \cdot \gamma_{hi}^{D,t_a}(k)}{\left(\sum_{\delta',\theta'} \gamma_{hi}^{\delta',\theta'}(k) \cdot N_h^{\delta',\theta'}(k) \right)^2 }  & \textrm{if } &  s_{hi}=1 \wedge \delta, \theta \neq D, t_a \\
            \Delta_h(N_h (k)) \cdot \frac{N_h^{D,t_a}(k)}{N_h (k)} & \textrm{if } & s_{hi}=0 \wedge \delta, \theta = D, t_a  \\             
            0 & \textrm{if } & s_{hi}=0 \wedge \delta, \theta \not= D, t_a  
        \end{array}
   \right.
\end{equation}

Similarly, we express  $q_{i\ell}^{D,t_a}$ as:
$$
  q_{i\ell}^{D,t_a} = \min \left( \beta_{i\ell} \cdot \Sigma_l(N_l) \cdot \frac{\gamma_{i\ell}^{D,t_a}\cdot N_i^{D,t_a}}{\sum_{\delta,\theta} \gamma_{i\ell}^{\delta,\theta} \cdot N_i^{\delta,\theta}}, \,
  \Delta_i(N_i) \cdot \frac{\gamma_{i\ell}^{D,t_a} \cdot N_i^{D,t_a}}{N_i} 
  \right)
$$
\noindent and as previously $s_{i\ell}=1$ when $q_{i\ell}^{D,t_a}$ is equal to the supply side in the above $\min$ and $s_{i\ell}=0$ when $q_{i\ell}^{D,t_a}$ is equal to the demand side in the above $\min$.

The partial derivatives of $q_{i\ell}^{D,t_a}$ concerning $\gamma_{i\ell}^{\delta,\theta}$ follows from (\ref{Eq_quvDta}):

\begin{equation}
\partial_{\gamma_{i\ell}^{\delta,\theta}(k)} q_{i\ell}^{D,t_a}(k) = 
\begin{cases}
\beta_{i\ell}\cdot \Sigma_\ell(N_\ell (k)) \cdot \frac{N_i^{D,t_a}(k)}{\sum_{\delta',\theta'} \gamma_{i\ell}^{\delta',\theta'}(k) \cdot N_i^{\delta',\theta'}(k)} \cdot \Bigg[1-\frac{N_i^{D,t_a}(k) \cdot \gamma_{i\ell}^{D,t_a}(k)}{\sum_{\delta',\theta'} \gamma_{i\ell}^{\delta',\theta'}(k) \cdot N_i^{\delta',\theta'}(k)} \Bigg] & \textrm{if } s_{hi}=1 \wedge \delta, \theta = D, t_a \\
-\beta_{i\ell}\cdot \Sigma_\ell(N_\ell (k))\cdot \frac{N_i^{D,t_a}(k) \cdot N_i^{\delta,\theta}(k) \cdot \gamma_{i\ell}^{D,t_a}(k)}{\Big(\sum_{\delta',\theta'} \gamma_{i\ell}^{\delta',\theta'}(k) \cdot N_i^{\delta',\theta'}(k) \Big)^2} & \textrm{if } s_{i\ell}=1 \wedge \delta, \theta \neq D, t_a \\
\Delta_i(N_i (k)) \cdot \frac{N_i^{D,t_a}}{N_i}(k) & \textrm{if } s_{i\ell}=0 \wedge \delta, \theta = D, t_a \\
0 & \textrm{if } s_{i\ell}=1 \wedge \delta, \theta = D, t_a
\end{cases} 
\end{equation}

\noindent Finally, we calculate the gradient of $E_i(k)$. Recall that:

\begin{equation}
E_i(k) = N_i(k) + \Bigg[\sum_{l \in \Gamma^-(i)} q_{li}(k)-\sum_{j \in \Gamma^+(i)} q_{ij}(k)\Bigg]\cdot \Delta t 
\end{equation}

\noindent and,

\begin{equation}
q_{uv}(k) = \min\Bigg[\beta_{uv} \cdot \Sigma_v(N_v(k)), \Delta_u(N_u(k))\cdot \sum_{\theta,\delta} \frac{\gamma_{uv}^{\delta,\theta}\cdot N_u^{\delta,\theta}(k)}{N_u(k)} \Bigg]    
\end{equation}
As previously $s_{uv}=1$ when $q_{uv}$ is equal to the supply side in the above $\min$ and $s_{uv}=0$ when $q_{uv}$ is equal to the demand side in the above $\min$.
The partial derivatives of $E_i(k)$ concerning $\gamma_{hi}^{\delta,\theta}$ and $\gamma_{i\ell}^{\delta,\theta}$ are:

\begin{equation}
\begin{cases}
\partial_{\gamma_{hi}^{\delta,\theta}(k)} E_i (k)= \Delta t \cdot \partial_{\gamma_{hi}^{\delta,\theta}(k)} q_{hi}(k) \\
\partial_{\gamma_{i\ell}^{\delta,\theta}(k)} E_i (k) = -\Delta t \cdot \partial_{\gamma_{i\ell}^{\delta,\theta}(k)} q_{i\ell}(k)
\end{cases}   
\end{equation}

\noindent where the partial derivatives of $q_{hi}$ and $q_{i\ell}$ are:

\begin{equation}
\partial_{\gamma_{hi}^{\delta,\theta}(k)} q_{hi}(k) =
\begin{cases}
0 & \textrm{if } s_{hi} = 1 \\
N_h^{\delta,\theta}(k) \cdot \frac{\Delta_h(N_h (k))}{N_h (k)}  & \textrm{if } s_{hi} = 0
\end{cases}
\end{equation}

\noindent and,

\begin{equation}
\partial_{\gamma_{i\ell}^{\delta,\theta}(k)} q_{i\ell}(k) =
\begin{cases}
0 & \textrm{if } s_{i\ell} = 1 \\
N_i^{\delta,\theta}(k) \cdot \frac{\Delta_i(N_i (k))}{N_i (k)}  & \textrm{if } s_{i\ell} = 0
\end{cases}
\end{equation}

This completes the calculation of the derivatives in step four.

Also note that if region $i$ is an Origin $O \in \mathcal{O}$, then $\partial_{\gamma_{hi}^{\delta,\theta}}=0$ because $\sum_{l \in \Gamma^-(i)} q_{li}(k)_{|i=0} $ is replaced by $\sum_{\delta,\theta} \Delta_O^{\delta,\theta}(k)$. 

\section{Projections}\label{Sect5}

\subsection{Projection on $\Cc$}\label{Sect5.1}

The domain $\Cc$ can be decomposed into a product of domains corresponding to the regions $i$, $\Cc = \times_{i,D,t_a} \Cc_i^{D,t_a}$:
\begin{equation}
\Cc_i^{D,t_a} \Def \begin{cases}
\sum_{j \in \Gamma^+(i)} \gamma_{ij}^{D,t_a} = 1 \\
\gamma_{ij}^{D,t_a} \ge 0, \forall j \in \Gamma^+(i)
\end{cases}
\end{equation}
The projection on $\Cc$ can be decomposed into projections specific to each region $i$, i.e. on the $\Cc_i^{D,t_a}$.

Let us calculate $\gamma = \mathcal{P}_{\Cc_i^{D,t_a}}(g)$. $\gamma$ is the solution of the following linear-quadratic problem

\begin{equation}\label{Eq_Proj0}
\min_{\gamma\in\Cc_i^{D,t_a}} \frac{1}{2}\cdot \sum_{j \in \Gamma^+(i)} |\gamma_{ij}^{D,t_a}-g_j|^2   
\end{equation}

As is usual in the case of simplex constraint domains, the KKT (Karush-Kuhn-Tucker) optimal conditions for  yielding the projection can be expressed as:

\begin{equation}\label{Eq_Proj1}
\gamma_{ij}^{D,t_a} = \mathcal{P}_{+}(g_j+\zeta)  
\end{equation}

\noindent with $\zeta$ being the KKT coefficient of the equality constraint in the projection $\Cc_i$ an $\mathcal{P}_{+}$ the projection on the set of positive numbers ($\mathcal{P}_{+}(y) = 0$ if $y\leq 0$ and $\mathcal{P}_{+}(y) = y$ if $y\geq 0$). The coefficient $\zeta$ is obtained by solving:

\begin{equation}\label{Eq_Proj2}
\sum_{j \in \Gamma^+(i)} \mathcal{P}_{+}(g_j+\zeta)=1     
\end{equation}

\begin{figure}[h!]\centering
    \includegraphics[width=0.3\textwidth]{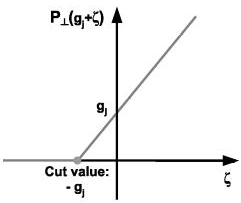}
    \caption{function $\zeta \rightarrow  \mathcal{P}_{+}(g_j+\zeta)$ with cut value $\zeta = -g_j$}
    \label{Fig:Projections_1}
\end{figure}

Figure \hyperref[Fig:Projections_1]{\ref*{Fig:Projections_1}} shows the function $\zeta \rightarrow  \mathcal{P}_{+}(g_j+\zeta)$. It illustrates the pointwise behaviour of the projection function $\mathcal{P}_{+}(g_j + \zeta)$ for a single component $j$ as a function of the KKT multiplier $\zeta$. The function exhibits a threshold structure, remaining at zero for $\zeta \le -g_j$ and increasing linearly with slope 1 beyond this cut value. This cut-off mechanism is central to the projection operation, as it governs which components of $g_j$ become active (positive) in the projected solution $\gamma$.

The function $\xi: \zeta \xrightarrow[]{\xi} \sum_{j \in \Gamma^+(i)} \mathcal{P}_{+}(g_j+\zeta)$ is convex increasing piecewise affine and invertible on $\xi^{-1}(\RR_{*}^+)$. Therefore, the solution of $\xi(\zeta)=1$ exists and is unique.

\begin{figure}[h!]\centering
    \includegraphics[width=0.4\textwidth]{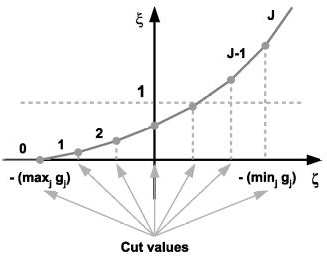}
    \caption{Function $\xi$ and its intersection with $1$.}
    \label{Fig:Projections_2}
\end{figure}

Figure \hyperref[Fig:Projections_2]{\ref*{Fig:Projections_2}} shows the shape of the function $\xi$. It depicts the aggregated function $\xi(\zeta) = \sum_{j \in \Gamma^+(i)} \mathcal{P}_{+}(g_j + \zeta)$, which must equal 1 to satisfy the normalisation constraint of the simplex domain $\Cc_i$. The function $\xi$ is piecewise affine, continuous, and convex, with breakpoints at each $-g_j$. The slope of $\xi$ increases in integer steps, corresponding to the number of active elements in the projection. The unique root $\zeta^*$ such that $\xi(\zeta^*) = 1$ ensures that the projected vector $\gamma$ lies on the probability simplex. The interpolation steps and structure illustrated in Figure~\ref{Fig:Projections_2} justify the algorithmic solution described above for computing the projection efficiently with complexity $\mathcal{O}(F \log F)$.

As Figure \hyperref[Fig:Projections_2]{\ref*{Fig:Projections_2}} shows, the function $\xi$ is continuous, convex, piecewise affine, and has an integer slope, which increases by one unit at each cut value $-g_j$. Thus, the slope of $xi$ is 0 if $\zeta \le -(\max_{j \in \Gamma^+(i)} g_j)$, and equal to the number of elements in $\Gamma^+(i)$ if $\zeta \ge -(\min_{j \in \Gamma^+(i)} g_j)$. ( $-(\max_{j \in \Gamma^+(i)} g_j)$ and $-(\min_{j \in \Gamma^+(i)} g_j)$ are the extreme cut values). 

The algorithm recommended for solving $\zeta$, i.e. $\xi(\zeta)=1$, runs as follows:

\begin{enumerate}
\item rank the cut values by increasing order: $-g_1 \le -g_2 \le \dots \le g_{\#\Gamma^+(i)}$;
\item find $\ell$ such that $\xi(-g_\ell)\le 1 < \xi(-g_{\ell+1})$;
\item interpolate:
\begin{itemize}
    \item Case 1: $\xi(-g_\ell)\le 1 < \xi(-g_{\ell+1}) $ then:
    \begin{equation}\label{Eq_Proj3.1}
    \zeta = -g_\ell + (g_\ell + g_{\ell+1})\cdot \frac{1-\xi(-g_\ell)}{\xi(-g_{\ell+1})-\xi(-g_\ell)}
    \end{equation}
As shown by Figure \hyperref[Fig:Projections_2]{\ref*{Fig:Projections_2}} we can carry out the following simplification in (\ref{Eq_Proj3.1}):
$$
    \frac{g_\ell + g_{\ell+1}}{\xi(-g_{\ell+1})-\xi(-g_\ell)} = \frac{1}{\ell}
$$
    \item Case 2: $\xi(-g_F) \le 1$, where $F$ denotes $\#\Gamma^+(i)$. If $\zeta \ge -g_F$, $\xi(\zeta) = \xi(-g_F) + J(\zeta+g_F)$. Thus:
    \begin{equation}\label{Eq_Proj3.2}
    \zeta = -g_F + (1+\xi(-g_F))/F    
    \end{equation}
\end{itemize}
\end{enumerate}

The complexity of the calculation of the solution of equation \hyperref[Eq_Proj2]{\ref*{Eq_Proj2}}, and thus the calculation of the projector on $\Cc_i$, via Equation \hyperref[Eq_Proj1]{\ref*{Eq_Proj1}}, is of the order of $F \ln F$, where $F = \#\Gamma^+(i)$.

\subsection{Projection on the domain $\Bc$}\label{Sect5.2}
The decision variables $\Delta (.)$ satisfy the total demand constraints \hyperref[Eq_TotalDemand1]{\ref*{Eq_TotalDemand1}} and positivity constraints, i.e. satisfy the set of constraints $(\Bc)$, which according to Equation \hyperref[Eq_Constraint3]{\ref*{Eq_Constraint3}} can by decomposed into a product od elementary convex simplex domains $\Bc^{D,t_a}_O$.
The domain $\Bc^{D,t_a}_O$, in a time-discretised setting, is given by:

\begin{equation}
\bar{\Delta}^{D,t_a}_O = \Bigg(\sum\limits_{k=0}^{K-1} \Delta^{D,t_a}_O(k) \Bigg) \cdot \frac{\Delta t}{|\mathcal{T}_D |}, \forall \Delta^{D,t_a}_O(k) \ge 0 \wedge \forall k=0,\dots,K-1
\end{equation}

Calculate of the projector $P_\Bc (\Delta) $. Formally, this calculation is similar to the calculation of $\mathcal{P}_{\mathcal{\Cc}}(.)$, because $\Bc$ is the product of independent sets $\Bc^{D,t_a}_O$ of simplex constraints (positivity constraint plus a single linear constraint $\bar{\Delta}^{D,t_a}_O$). Thus, $\Delta^{D,t_a}_O = P_{\Bc^{D,t_a}_O} \left( d_O^{D,t_a} \right) $ is given by:

\begin{equation}
\Delta^{D,t_a}_O(k) = \mathcal{P}_{+} (d_O^{D,t_a}(k)+\zeta)    
\end{equation}

\begin{equation}
\bar{\Delta}^{D,t_a}_O \cdot \frac{|\tau_D|}{\Delta t} = \sum\limits_{k=0}^{K-1} \Delta_O^{D,t_a}(k) 
\end{equation}

The term $\zeta$ is a KKT coefficient solved as

\begin{equation}\label{Eq_Proj5}
\bar{\Delta}^{D,t_a}_O \cdot \frac{|\tau_D|}{\Delta t} = \sum\limits_{k=0}^{K-1} \mathcal{P}_{+}(d_O^{D,t_a}(k)+\zeta)
\end{equation}

\noindent which is formally identical to Equation \hyperref[Eq_Proj2]{\ref*{Eq_Proj2}}. 

We define the function $\eta$ of $\zeta$ as follows:

\begin{equation}
\eta := \zeta \xrightarrow[]{\eta} \sum\limits_{k=0}^{K-1} \mathcal{P}_{+}(d_O^{D,t_a}(k)+\zeta) 
\end{equation}

\noindent The function $\eta$ is convex, increasing, piecewise, and affine. The function $\zeta \rightarrow \mathcal{P}_{+}(d_O^{D,t_a}(k)+\zeta)$ admits $-d_O^{D,t_a}(k)$ as a cut value (value of $\zeta$ for which this function changes slope). Order the cut values $-d_O^{D,t_a}(k)$ from the smallest to the largest values. We denote $\bar{\omega}$ as the corresponding permutation of $\{0,1,\dots,K-1\}$:

\begin{equation}\label{Eq_Proj6}
-d_O^{D,t_a}(\bar{\omega}(0)) \le -d_O^{D,t_a}(\bar{\omega}(1)) \le \dots \le -d_O^{D,t_a}(\bar{\omega}(K-1))   
\end{equation}

\begin{figure}[h!]\centering
    \includegraphics[width=0.4\textwidth]{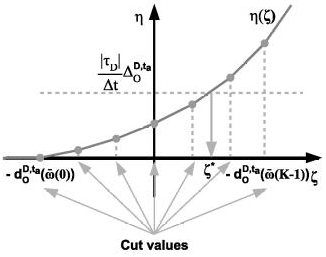}
    \caption{Intersection of $\bar{\Delta}^{D,t_a}_O \cdot \frac{|\tau_D|}{\Delta t} $ with $\eta(\zeta)$}
    \label{Fig:Projections_3}
\end{figure}

Figure \hyperref[Fig:Projections_3]{\ref*{Fig:Projections_3}} shows the function $\eta$ and its intersection with $\bar{\Delta}^{D,t_a}_O \cdot \frac{|\tau_D|}{\Delta t} $.

The recommended algorithm for solving equation \hyperref[Eq_Proj5]{\ref*{Eq_Proj5}} is similar to the algorithm for $\mathcal{P}_{\alpha_j}$:

\begin{enumerate}
    \item rank the cut values in increasing order: find $\bar{\omega}$ permutations of $\{0,1,\dots,K-1\}$ that satisfy Equation \hyperref[Eq_Proj6]{\ref*{Eq_Proj6}}.
    \item find $l$ such that $\eta(-d_O^{D,t_a}(\bar{\omega}(l))) \le \bar{\Delta}^{D,t_a}_O \cdot \frac{|\tau_D|}{\Delta t} \le \eta(-d_O^{D,t_a}(\bar{\omega}(l+1)))$.
    \item interpolate:
    \begin{itemize}
        \item Case 1, when $\eta(-d_O^{D,t_a}(\bar{\omega}(l))) \le \bar{\Delta}^{D,t_a}_O \cdot \frac{|\tau_D|}{\Delta t} \le \eta(-d_O^{D,t_a}(\bar{\omega}(l+1)))$:
        \begin{equation}
            \zeta = -d_O^{D,t_a}(\bar{\omega}(l)) + \Big(1-\eta(-d_O^{D,t_a}(\bar{\omega}(l)))\Big) \cdot \frac{d_O^{D,t_a}(\bar{\omega}(l))-d_O^{D,t_a}(\bar{\omega}(l+1))}{\eta(-d_O^{D,t_a}(\bar{\omega}(l+1)))-\eta(-d_O^{D,t_a}(\bar{\omega}(l)))}
        \end{equation}
        \item Case 2, when $\bar{\Delta}^{D,t_a}_O \cdot \frac{|\tau_D|}{\Delta t} \ge \eta(-d_O^{D,t_a}(\bar{\omega}(K-1)))$:
        \begin{equation}
            \zeta = -d_O^{D,t_a}(\bar{\omega}(K-1)) + [1-\eta(d_O^{D,t_a}(\bar{\omega}(K-1)))]\cdot K
        \end{equation}
        \item Case 3, when $-d_O^{D,t_a}(\bar{\omega}(l)) = -d_O^{D,t_a}(\bar{\omega}(l+1))$:
        \begin{equation}
            \zeta = -d_O^{D,t_a}(\bar{\omega}(l)) = -d_O^{D,t_a}(\bar{\omega}(l+1))
        \end{equation}
    \end{itemize}
\end{enumerate}

Note that the third case must be added since the sequence of cut values may not satisfy the monotonic condition after the reordering.

\section{Results}\label{Sect6}

To evaluate the effectiveness of the proposed SO control framework, we compare its performance against the state-of-the-art gradient-based method that uses pairwise marginal travel costs \citep{Yildirimoglu-2018, ameli2020improving}. For benchmarking purposes, we also include two established solution algorithms: the MSA \citep{mahmassani1993network} and a gap-based approach \citep{lu2009equivalent}. Numerical experiments are conducted on a synthetic network representing a large-scale abstraction of the Paris metropolitan area, divided into eight interconnected regions. Each region is modelled using a triangular Macroscopic Fundamental Diagram (MFD), defined by a piecewise linear speed-density relationship.

\begin{figure}[!h]
\centering
\begin{tikzpicture}[scale=1.2, every node/.style={circle, draw, minimum size=1cm, font=\small, fill=white}]

\node (R1) at (0, 2.5) {R1};
\node (R2) at (2.5, 0.5) {R2};
\node (R3) at (5, 0.5) {R3};
\node (R4) at (1.75, 2.5) {R4};
\node (R5) at (3.75, 2.5) {R5};
\node (R6) at (5.5, 2.5) {R6};
\node (R7) at (2.5, 4.5) {R7};
\node (R8) at (4.5, 4.5) {R8};

\draw[-] (R1) -- (R2);
\draw[-] (R1) -- (R4);
\draw[-] (R1) -- (R7);
\draw[-] (R2) -- (R3);
\draw[-] (R2) -- (R4);
\draw[-] (R2) -- (R5);
\draw[-] (R3) -- (R6);
\draw[-] (R4) -- (R5);
\draw[-] (R4) -- (R7);
\draw[-] (R5) -- (R6);
\draw[-] (R5) -- (R7);
\draw[-] (R5) -- (R3);
\draw[-] (R5) -- (R8);
\draw[-] (R6) -- (R8);
\draw[-] (R7) -- (R8);

\end{tikzpicture}
\caption{Schematic layout of the 8-region network used in the experiments.}
\label{fig:region-network}
\end{figure}
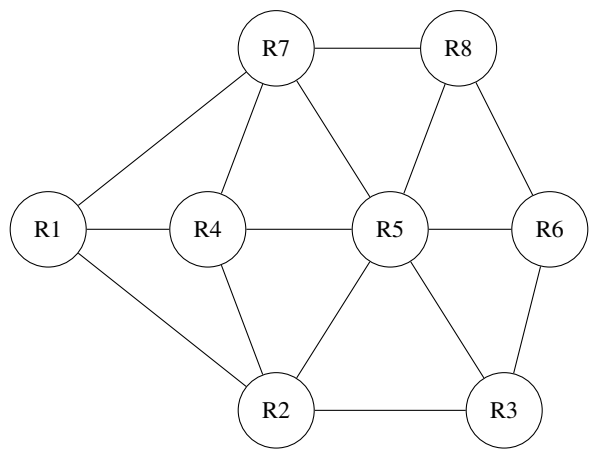

Figure~\ref{fig:region-network} depicts the network topology, where nodes represent urban regions governed by MFDs, and edges represent interregional travel links. The central region, R5, represents the core of Paris within the Île-de-France area and serves as a major origin and destination hub for the morning peak-hour scenario. Table~\ref{tab:parameters} lists the simulation parameters. The simulations span a four-hour period (06:00–10:00 am), featuring a trapezoidal inflow profile with peak demand occurring between 07:00 and 09:00. Although the actual region sees approximately 1.5 million weekday morning car trips, the total simulated demand is rescaled to 150,000 vehicle trips for computational feasibility. All MFD parameters are adjusted accordingly. Trip lengths range from 2 to 8 regions, with desired arrival times uniformly distributed between 08:00 and 10:15.

The origin-destination (OD) demand distribution is shown in Table~\ref{tab:od-matrix}. The majority of trips either occur within R5 or originate in neighbouring suburban regions (R1–R4 and R6–R8) with destinations in the Paris core. Internal trips within R5 account for 15\% of the total demand, reflecting typical intra-city travel. Suburban regions such as R1, R2, and R3 are major contributors to inbound flows into R5, consistent with observed morning commuting patterns.

\begin{table}[!h]
\centering
\caption{Simulation Parameters}
\label{tab:parameters}
\begin{tabular}{@{}ll@{}}
\toprule
\textbf{Parameter}               & \textbf{Value}             \\ \midrule
Number of regions               & 8                          \\
Number of vehicle trips         & 150000                       \\
Simulation time (Tsimul)        & 14400 seconds (6:00 am to 10:00 am)      \\
Time step (Dt)                  & 1 second                   \\
Free flow speed (R1-R3 and R7-R8)                 & 24 m/s   \\
Free flow speed (R4 and R6)                       & 14 m/s   \\
Free flow speed (R1-R3 and R7-R8)                 & 10 m/s   \\
Critical accumulation ($n_c$)   & 3000 vehicles              \\
Jam density ($n_j$)             & 12000 vehicles              \\
Segment Max length (R1-R4 and R6-R8)             & 15 km      \\
Segment Max length (R5)                          & 10 km      \\
Demand inflow profile           & Trapezoidal shape (peak at 7:00 am to 9:00 am) \\
Desired arrival time window      & 10–20 minutes              \\
Route length range              & 2–8 regions per trip       \\ \bottomrule
\end{tabular}
\end{table}

\begin{table}[!h]
\centering
\caption{Origin-Destination (OD) matrix of regional demand. Each cell shows the percentage of total trips from an origin region (row) to a destination region (column). Desired arrival times are uniformly distributed between 8:00 AM and 9:30 AM.}
\label{tab:od-matrix}
\begin{tabular}{lcccccccc}
\toprule
\textbf{O / D} & \textbf{R1} & \textbf{R2} & \textbf{R3} & \textbf{R4} & \textbf{R5} & \textbf{R6} & \textbf{R7} & \textbf{R8} \\ \midrule
\textbf{R1} & 2\% & 0\% & 0\% & 0\% & 10\% & 0\% & 0\% & 0\% \\
\textbf{R2} & 0\% & 2\% & 0\% & 0\% & 9\% & 0\% & 0\% & 0\% \\
\textbf{R3} & 0\% & 0\% & 2\% & 0\% & 8\% & 0\% & 0\% & 0\% \\
\textbf{R4} & 0\% & 0\% & 0\% & 2\% & 7\% & 0\% & 0\% & 0\% \\
\textbf{R5} & 4\% & 3\% & 0\% & 0\% & \textbf{15\%} & 3\% & 2\% & 0\% \\
\textbf{R6} & 0\% & 0\% & 0\% & 0\% & 6\% & 2\% & 0\% & 0\% \\
\textbf{R7}& 0\% & 0\% & 0\% & 0\% & 6\% & 0\% & 2\% & 0\% \\
\textbf{R8} & 0\% & 0\% & 0\% & 0\% & 5\% & 0\% & 0\% & 2\% \\ \bottomrule
\end{tabular}
\end{table}

The proposed $\beta$-controller algorithm was applied to the test network and benchmarked against two established methods: the baseline MSA and the gap-based algorithm. Detailed formulations of MSA and the gap-based approach are provided in \cite{ameli2021computational}. The SO $\beta$-controller seeks to minimise total system travel time by optimising departure rate allocations through iterative estimation of marginal costs across all feasible paths. In contrast, both MSA and the gap-based method estimate path marginal costs via the aggregation of link-level marginal costs. All methods initialise from a uniform departure time distribution relative to the desired arrival time window.

Table~\ref{tab:traffic_assignment_comparison} provides a detailed comparison of final results for the SO $\beta$-controller, MSA, and Gap Method. The SO $\beta$-controller outperforms both benchmark methods, reducing total travel cost by 14\% relative to MSA and by 12\% compared to the gap-based method. Region-specific improvements are also evident. In region R5 (Paris), the maximum vehicle accumulation is reduced to 9419 under the SO $\beta$-controller, compared to 9785 with MSA and 9276 with the gap method. This reduction is accompanied by an increase in average speed: 6.37 m/s under the SO controller, representing a 14\% improvement over MSA (5.60 m/s) and a 25\% gain over the gap-based method (5.11 m/s). Similar benefits are observed in other highly congested regions such as R4 and R6, where the SO controller reduces peak accumulations by more than 11\% and increases average speeds by 7–12\%.

In less congested regions (R1, R2, and R3), the SO $\beta$-controller maintains higher average speeds (exceeding 15.8 m/s) while keeping accumulations consistently below those of the benchmark methods. These results highlight the controller’s capacity to proactively manage demand, mitigate queuing in critical areas, and achieve a more balanced network flow than traditional equilibrium-based approaches.

\begin{table}[htbp]
\centering
\caption{Summary of final results for each algorithm: Avg Cost: average traveller cost; MxAcumm: maximum accumulation (veh); ASpeed: average speed (m/s).}
\begin{adjustbox}{width=\textwidth}
\begin{tabular}{lcc*{8}{cc}}
\hline
\multirow{2}{*}{\textbf{Algorithm}} & \multirow{2}{*}{\textbf{Total Travel Cost}} & \multirow{2}{*}{\textbf{Avg Cost (euros)}} 
& \multicolumn{2}{c}{\textbf{R1}} & \multicolumn{2}{c}{\textbf{R2}} & \multicolumn{2}{c}{\textbf{R3}} & \multicolumn{2}{c}{\textbf{R4}} \\
\cline{4-11}
 & & & \textbf{MxAccum} & \textbf{ASpeed} & \textbf{MxAccum} & \textbf{ASpeed} & \textbf{MxAccum} & \textbf{ASpeed} & \textbf{MxAccum} & \textbf{ASpeed} \\
\hline
\textbf{SO $\beta$-controller} & 536,908 & 3.58 & 5580 & 15.82 & 5153 & 16.52 & 4793 & 17.44 & 5658 & 14.89 \\
\textbf{MSA} & 624,333 & 4.16 & 6068 & 14.78 & 5595 & 15.34 & 5198 & 16.71 & 6416 & 13.06 \\
\textbf{Gap Method} & 611,752 & 4.08 & 6571 & 13.49 & 5888 & 14.97 & 5506 & 15.23 & 6637 & 11.61 \\
\hline
\multicolumn{3}{c}{} & \multicolumn{2}{c}{\textbf{R5}} & \multicolumn{2}{c}{\textbf{R6}} & \multicolumn{2}{c}{\textbf{R7}} & \multicolumn{2}{c}{\textbf{R8}} \\
\cline{4-11}
\textbf{Algorithm} & \textbf{Total Travel Cost} & \textbf{Avg Cost (euros)} & \textbf{MxAccum} & \textbf{ASpeed} & \textbf{MxAccum} & \textbf{ASpeed} & \textbf{MxAccum} & \textbf{ASpeed} & \textbf{MxAccum} & \textbf{ASpeed} \\
\hline
\textbf{SO $\beta$-controller} & 536,908 & 3.58 & 9419 & 6.37 & 5437 & 14.03 & 5066 & 18.45 & 5313 & 16.33 \\
\textbf{MSA} & 624,333 & 4.16 & 9785 & 5.60 & 5856 & 12.63 & 5527 & 16.43 & 5431 & 15.74 \\
\textbf{Gap Method} & 611,752 & 4.08 & 9276 & 5.11 & 6275 & 12.32 & 5734 & 15.84 & 5843 & 14.64 \\
\hline
\end{tabular}
\end{adjustbox}
\label{tab:traffic_assignment_comparison}
\end{table}

\begin{figure}[htbp]
    \centering
    \includegraphics[width=0.85\textwidth]{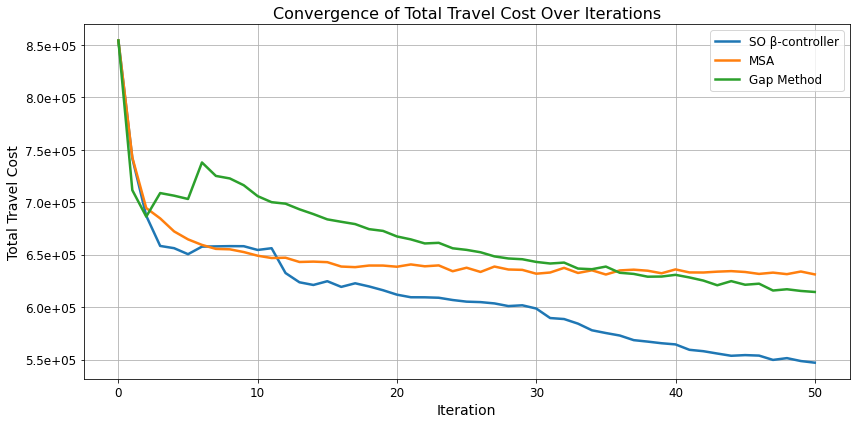}
    \caption{Convergence of total travel cost over 50 iterations for the three traffic assignment methods: SO $\beta$-controller, MSA, and Gap Method. All methods begin from a common initial cost. The SO $\beta$-controller achieves the fastest and most stable convergence, reaching the lowest final travel cost. MSA displays a steady, gradual decline as it systematically adjusts non-optimal departure times. The Gap Method shows higher initial oscillations, attributable to its trip, based update mechanism—but ultimately stabilises at a lower cost than MSA, though still above that of the SO $\beta$-controller.}
    \label{fig:convergence_plot}
\end{figure}

Figure~\ref{fig:convergence_plot} illustrates the convergence behaviour of total travel cost over iterations for the three traffic assignment methods. The SO $\beta$-controller exhibits the fastest and most stable convergence, achieving a substantial reduction in total cost within the first 10 iterations and stabilising around 537,000 by iteration 50. In comparison, the MSA converges more gradually, reaching a steady-state value near 631,000. The gap-based method initially shows a sharp decline in cost but is characterised by persistent oscillations, ultimately converging to a higher cost level of approximately 614,000. These results underscore the superior convergence efficiency of the SO $\beta$-controller, which more effectively leverages marginal cost feedback to achieve faster and more optimal coordination across the network.

\section{Conclusions}\label{Sect7}

This paper introduced a SO-based traffic control framework based on a $\beta$-controller, designed to dynamically allocate departure rates to minimise total travel cost across a multi-region urban network governed by the MFD dynamics. The proposed approach was benchmarked against two established traffic assignment algorithms - the MSA and a gap-based method—on a synthetic eight-region network representative of the Paris Île-de-France area. The testbed incorporated realistic spatial and temporal demand distributions, as well as peak-period congestion dynamics.


The core innovation of the framework lies in its integration of a path-based marginal cost control mechanism that directly leverages MFD feedback. This enables faster convergence and more effective network coordination compared to traditional methods. Numerical experiments showed that the SO $\beta$-controller consistently outperformed both benchmark algorithms, achieving reductions of up to 14\% in total system travel cost, improving average user costs, and mitigating congestion, particularly in critical regions such as central Paris (R5), where it maintained higher travel speeds and reduced maximum vehicle accumulations.


Despite its promising performance, the current framework has several limitations. First, it assumes deterministic demand and fixed route sets, which may not fully capture the stochastic and adaptive nature of real-world traveller behaviour. Second, the model relies on pre-defined MFDs and assumes homogeneous compliance with assigned departure schedules, potentially limiting its applicability in heterogeneous or partially regulated networks. Finally, the current formulation focuses exclusively on departure time optimisation, without accounting for elastic demand or pricing mechanisms.


Future research will aim to address these limitations by extending the framework to incorporate demand uncertainty, adaptive routing behaviour, and network disruption scenarios. Integrating time-varying MFDs and dynamic pricing strategies may further enhance the model’s realism and policy relevance. In addition, future work will involve applying the $\beta$-controller to real-world case studies, supported by empirical calibration of MFDs, to validate its effectiveness for metropolitan-scale traffic management.


\section*{Acknowledgments}
\begin{sloppypar}

M. Ameli acknowledges support from the French ANR research project SMART-ROUTE (grant number ANR-24-CE22-7264).

M. Men\'endez acknowledges support from the NYUAD Center for Interacting Urban Networks (CITIES) funded by Tamkeen under the NYUAD Research Institute Award CG001.
\end{sloppypar}

\appendix
\section{Alternative origin model}\label{Sect6.3}
The idea is to avoid considering the origin cells. Therefore, we keep the following.

\begin{itemize}
    \item state variables: $N_O(k)$ and $N_O^{D,t_a}$ ($\forall k, \forall d \in \mathcal{D}, \forall t_a \in \mathcal{T}_a$);
    \item commands: $\gamma_{Oi}^{D,t_a}$ and $\Delta_O^{D,t_a}(k)$ ($\forall i \in \Gamma^+(O), \forall d \in \mathcal{D}, \forall t_a \in \mathcal{T}_a$)
\end{itemize}

We treat $N_O$ as a homogeneous queue:

\begin{enumerate}
    \item Demand:
    \begin{equation}\label{Eq_M4}
        \delta_O(k) = \Delta_O(k) + \frac{N_O(k)}{\Delta t}
    \end{equation}
    \item Partial demands:
    \begin{equation}\label{Eq_M51}
        \begin{cases}
            \gamma_{Oi}(k) \stackrel{\text{def}}{=} \frac{\sum_{\substack{\delta \in {\mathcal{D}} \\ \theta \in \mathcal{T}_a}} \gamma_{Oi}^{\delta,\theta}\cdot \left(\Delta_O^{\delta,\theta}(k)\cdot \Delta t + N_O^{\delta,\theta}(k)\right)}{\Delta_O(k)\cdot \Delta t + N_O(k)} \\
            \delta_{Oi}(k) = \gamma_{Oi}(k) \cdot \delta_O(k)
        \end{cases}
    \end{equation}
\end{enumerate}

Note that the second equation in the previous system of equations \hyperref[Eq_M51]{\ref*{Eq_M51}} can be expressed as:

\begin{equation}
\delta_{Oi}(k) = \sum_{\substack{\delta \in \mathcal{D} \\ \theta \in \mathcal{T}_a }} \gamma_{Oi}^{\delta,\theta}\cdot \Bigg(\Delta_O^{\delta,\theta}(k) + \frac{N_O^{\delta,\theta}(k)}{\Delta t}\Bigg), \forall i \in \Gamma^+(O)
\end{equation}

Then, we calculate the outflows as follows:

\begin{equation}\label{Eq_M6}
q_{Oi}(k) = \min(\delta_{Oi}(k),\sigma_{Oi}(k))    
\end{equation}

\noindent and then the conservation equation for the queue $N_O$ is:

\begin{equation}\label{Eq_M7}
N_O(k+1) = N_O(k) + \Delta_O(k)\cdot \Delta t - \sum\limits_{i \in \Gamma^+(O)} q_{Oi}(k)\cdot \Delta t = E_O(k)    
\end{equation}

\begin{equation}\label{Eq_M8}
    \begin{cases}
        N_O^{D,t_a}(k+1) = N_O^{D,t_a}(k) + \Delta t \cdot (\Delta_O^{D,t_a}(k) - \sum\limits_{i \in \Gamma^+(O)} q_{Oi}(k))\\
        q_{Oi}^{D,t_a}(k) = q_{Oi}(k) \cdot \frac{\Big(\Delta_O^{D,t_a}(k)\cdot \Delta t + N_O^{D,t_a}(k)\Big)\cdot \gamma_{Oi}^{D,t_a}(k)}{\sum_{\substack{\delta \in \mathcal{D} \\ \theta \in \mathcal{T}_a}} \Big(\Delta_O^{\delta,\theta}(k)\cdot \Delta t + N_O^{\delta,\theta}(k) \Big)\cdot \gamma_{Oi}^{\delta,\theta}(k)} = E_O^{D,t_a}(k)
    \end{cases}
\end{equation}

In this scheme, there is also a significant amount of smoothing, especially in the system of equations defined in \hyperref[Eq_M8]{\ref*{Eq_M8}} (no track of the order in which travellers enter the queue for $D$ and $t_a$ is maintained). The term $\gamma_{Oi}^{D,t_a}$ is optimised for the global composition of the queue $N_O$. If this queue is large, this implies that there is a strong smoothing effect. On the other hand, if $N_O$ is small, then the scheme defined in Equations \hyperref[Eq_M4]{\ref*{Eq_M4}}, \hyperref[Eq_M51]{\ref*{Eq_M51}}, \hyperref[Eq_M6]{\ref*{Eq_M6}}, \hyperref[Eq_M7]{\ref*{Eq_M7}} and \hyperref[Eq_M8]{\ref*{Eq_M8}}, is probably much better (no congestion at the entrance of the network), because then there is no smoothing.

\bibliography{references}

\nolinenumbers

\end{sloppypar}

\end{document}